\documentstyle[12pt]{article}
\def\Bbb R{{\rm \bf R}}
\def\proclaim#1{\vskip2mm{\bf #1}\em}
\def\endproclaim{\em \vskip2mm}
\def\tag#1{\eqno(#1)}
\def\gathered{\begin{array}{c}}
\def\endgathered{\end{array}}
\def\text{\mbox}

\begin{document}

\title {On finding an obstacle with the Leontovich boundary condition
via the time domain enclosure method}
\author{Masaru IKEHATA\footnote{
Laboratory of Mathematics,
Institute of Engineering,
Hiroshima University,
Higashi-Hiroshima 739-8527, JAPAN}}
%\date{}
\maketitle

\begin{abstract}
An inverse obstacle scattering problem for the wave governed by the Maxwell system
in the time domain, in particular, over a finite time interval  is considered.
It is assumed that the electric field $\mbox{\boldmath $E$}$ and magnetic field
$\mbox{\boldmath $H$}$ which are solutions of the Maxwell system
are generated only by a current density at the initial time located not far a way from an unknown obstacle.
The obstacle is embedded in a medium like air which has constant electric permittivity $\epsilon$
and magnetic permeability $\mu$.
It is assumed that the fields on the surface of the obstacle satisfy
the impedance-or the Leontovich boundary condition $\mbox{\boldmath $\nu$}\times\mbox{\boldmath $H$}
-\lambda\,\mbox{\boldmath $\nu$}\times(\mbox{\boldmath $E$}\times\mbox{\boldmath $\nu$})=\mbox{\boldmath $0$}$
with $\lambda$ an unknown positive function and $\mbox{\boldmath $\nu$}$
the unit outward normal.
The observation data are given by the electric field observed at the same place as the support of
the current density over a finite time interval.  It is shown that
an indicator function computed from the electric fields corresponding two current densities
enables us to know: the distance of the center of the common spherical support of the current densities to the obstacle; whether the value of the impedance $\lambda$
is greater or less than the special value $\sqrt{\epsilon/\mu}$.

\noindent
AMS: 35R30, 35L50, 35Q61, 78A46, 78M35

\end{abstract}

%\tableofcontents

\section{Introduction}

In this paper, we consider an inverse obstacle scattering problem for the wave
governed by the Maxwell system in the time domain, in particular, over a finite time interval.
The formulation of the problem basically follows that of \cite{IEM}.

We assume that the electric field $\mbox{\boldmath $E$}$ and magnetic field $\mbox{\boldmath $H$}$
are generated only by the current density $\mbox{\boldmath $J$}$ at the initial time
located not far a way from an unknown obstacle.  The obstacle is embedded in
a medium like air which has constant electric permittivity $\epsilon(\,>0)$ and magnetic permeability $\mu(>0)$.
On the surface of the obstacle unlike \cite{IEM} it is assumed that
the fields $\mbox{\boldmath $E$}$ and $\mbox{\boldmath $H$}$ satisfy the impedance- or
the Leontovich boundary condition (\cite{CK, KH}).

Let us formulate the problem more precisely.
We denote by $D$ the unknown obstacle.
We assume that $D$ is a non empty bounded open set of $\Bbb R^3$ with $C^2$-boundary such that $\Bbb R^3\setminus\overline D$
is connected.

Let $0<T<\infty$.
The governing equations of $\mbox{\boldmath $E$}$ and $\mbox{\boldmath $H$}$ over the time interval
$]0,\,T[$
take the form
$$\displaystyle
\left\{
\begin{array}{ll}
\displaystyle
\epsilon\frac{\partial\mbox{\boldmath $E$}}{\partial t}
-\nabla\times\mbox{\boldmath $H$}=\mbox{\boldmath $J$}
& \text{in}\,(\Bbb R^3\setminus\overline D)\times\,]0,\,T[,\\
\\
\displaystyle
\mu\frac{\partial\mbox{\boldmath $H$}}{\partial t}
+\nabla\times\mbox{\boldmath $E$}=\mbox{\boldmath $0$}
& \text{in}\,(\Bbb R^3\setminus\overline D)\times\,]0,\,T[,
\\
\\
\displaystyle
\mbox{\boldmath $E$}\vert_{t=0}=\mbox{\boldmath $0$}
& \text{in}\,\Bbb R^3\setminus\overline D,
\\
\\
\displaystyle
\mbox{\boldmath $H$}\vert_{t=0}=\mbox{\boldmath $0$}
& \text{in}\,\Bbb R^3\setminus\overline D
\end{array}
\right.
\tag {1.1}
$$
and
$$
\displaystyle
\begin{array}{ll}
\mbox{\boldmath $\nu$}\times\mbox{\boldmath $H$}
-\lambda\,
\mbox{\boldmath $\nu$}\times(\mbox{\boldmath $E$}\times\mbox{\boldmath $\nu$})=\mbox{\boldmath $0$}
& \text{on}\,\partial D\times\,]0,\,T[.
\end{array}
\tag {1.2}
$$
Note that $\mbox{\boldmath $\nu$}$ denotes the unit outward normal to $\partial D$.
To ensure the solvability of the initial boundary value problem
(1.1) and (1.2) by the theory of $C_0$ contraction semigroup \cite{Y}
we assume that $\lambda\in C^1(\partial D)$ and satisfies $\lambda (x)\ge 0$ for all $x\in\partial D$.
Under the condition $\mbox{\boldmath $J$}\in
C^1([0,\,T],\,L^2(\Bbb R^3\setminus\overline D)^3)$
we have the unique solution $(\mbox{\boldmath $E$}, \mbox{\boldmath $H$})$
which belongs to $C^1([0,\,T],\,L^2(\Bbb R^3\setminus\overline D)^3\times L^2(\Bbb R^3\setminus\overline D)^3)$
with $(\nabla\times\mbox{\boldmath $E$}(t),\,\nabla\times\mbox{\boldmath $H$}(t))\in
L^2(\Bbb R^3\setminus\overline D)^3\times L^2(\Bbb R^3\setminus\overline D)^3$
and (1.2) is satisfied in the sense of the trace \cite{DL3}.

The boundary condition in (1.2) is called the Leontovich boundary condition
and equivalent to the condition
$$\begin{array}{ll}
\displaystyle
\mbox{\boldmath $\nu$}\times(\mbox{\boldmath $H$}\times\mbox{\boldmath $\nu$})+\lambda\,\mbox{\boldmath $\nu$}\times
\mbox{\boldmath $E$}
=\mbox{\boldmath $0$}
& \text{on}\,\partial D\times\,]0,\,T[.
\end{array}
$$
In what follows we use these equivalent forms without mentioning explicitly.

The mathematical role of the existence of the impedance $\lambda$
can be seen by formally differentiating the energy
$$\displaystyle
{\cal E}(t)
=\frac{1}{2}
\int_{\Bbb R^3\setminus\overline D}(\epsilon\vert\mbox{\boldmath $E$}\vert^2+\mu\vert\mbox{\boldmath $H$}\vert^2)dx.
$$
By the help of the divergence theorem and (1.2), we have
$$\displaystyle
{\cal E}'(t)
=-\int_{\partial D}\lambda\vert\mbox{\boldmath $\nu$}\times(\mbox{\boldmath $E$}\times\mbox{\boldmath $\nu$})\vert^2 dS
+\int_{\Bbb R^3\setminus\overline D}\mbox{\boldmath $J$}\cdot\mbox{\boldmath $E$}dx.
$$
Thus if $\mbox{\boldmath $J$}(t)=0$ for $t\in[\delta,\,T]$ with a small $\delta>0$, we have
${\cal E}'(t)\le 0$ therein.  The solution loses its energy on the surface of the obstacle.

There should be several choices of current density
$\mbox{\boldmath $J$}$ as a model of antenna \cite{B, CB}. In this
paper, as considered in \cite{IEM} we assume that $\mbox{\boldmath
$J$}$ takes the form
$$
\displaystyle
\mbox{\boldmath $J$}(x,t)=f(t)\chi_B(x)\mbox{\boldmath $a$},
\tag {1.3}
$$
where $\mbox{\boldmath $a$}$ is an arbitrary unit vector; $B$ is a (very small)
open ball satisfying $\overline B\cap\overline D=\emptyset$
and $\chi_B$ denotes the characteristic function of $B$; $f\in C^1[0,\,T]$ with $f(0)=0$.

We consider the following problem.

{\bf\noindent Problem.}  Fix a large (to be determined later) $T<\infty$.  Generate the solutions $\mbox{\boldmath $E$}$ and $\mbox{\boldmath $H$}$
of system (1.1) and (1.2) by the source $\mbox{\boldmath $J$}$ having the form (1.3) with $f\not=0$ and observe $\mbox{\boldmath $E$}$ on $B$
over the time interval $]0,\,T[$.  The unit vector $\mbox{\boldmath $a$}$ on (1.3)
should be taken from a set of three linearly independent vectors.
Extract information about the geometry of $D$ and the qualitative state of the distribution
of $\lambda$ over $\partial D$ from the observed data.

As the author knows there is no result for this problem. Note that
in \cite{IEM} the perfect conductive boundary condition
$\mbox{\boldmath $\nu$}\times\mbox{\boldmath $E$}=\mbox{\boldmath
$0$}$ on $\partial D$ which is the extreme case  $\lambda=+\infty$
of the Leontovich boundary condition has been considered. Thus,
therein extracting the geometry of $D$ is the main interest.  Here
we wish to know the {\it qualitative state} of the surface of the
obstacle which is described by the unknown function  $\lambda$ on
$\partial D$.

\subsection{The statement of the results}

Let $\mbox{\boldmath $E$}_0$ and $\mbox{\boldmath $H$}_0$ be the solutions of
$$\displaystyle
\left\{
\begin{array}{ll}
\displaystyle
\epsilon\frac{\partial\mbox{\boldmath $E$}}{\partial t}
-\nabla\times\mbox{\boldmath $H$}=\mbox{\boldmath $J$}
&  \text{in}\,\Bbb R^3\times\,]0,\,T[,
\\
\\
\displaystyle
\mu\frac{\partial\mbox{\boldmath $H$}}{\partial t}
+\nabla\times\mbox{\boldmath $E$}=\mbox{\boldmath $0$}
&  \text{in}\,\Bbb R^3\times\,]0,\,T[,
\\
\\
\displaystyle
\mbox{\boldmath $E$}\vert_{t=0}=\mbox{\boldmath $0$}
&  \text{in}\,\Bbb R^3,
\\
\\
\displaystyle
\mbox{\boldmath $H$}\vert_{t=0}=\mbox{\boldmath $0$}
&  \text{in}\,\Bbb R^3.
\end{array}
\right.
\tag {1.4}
$$

Using $\mbox{\boldmath $E$}$ and $\mbox{\boldmath $E$}_0$ on $B$ over time interval $]0,\,T[$,
we introduce the indicator function of the enclosure method in this paper
$$
\displaystyle
I_{\mbox{\boldmath $f$}}(\tau,T)
=\int_B\mbox{\boldmath $f$}(x,\tau)\cdot(\mbox{\boldmath $W$}_e-\mbox{\boldmath $V$}_e)dx,
\tag {1.5}
$$
where
$$
\displaystyle
\mbox{\boldmath $f$}(x,\tau)
=-\frac{\tau}{\epsilon}\int_0^Te^{-\tau t}\mbox{\boldmath $J$}(x,t)dt,
\tag {1.6}
$$
$$
\displaystyle
\mbox{\boldmath $W$}_e=\mbox{\boldmath $W$}_e(x,\tau)
=\int_0^T e^{-\tau t}\mbox{\boldmath $E$}(x,t)dt
\tag {1.7}
$$
and
$$
\displaystyle
\mbox{\boldmath $V$}_e=\mbox{\boldmath $V$}_e(x,\tau)
=\int_0^Te^{-\tau t}\mbox{\boldmath $E$}_0(x,t)dt.
\tag {1.8}
$$
Note that we have
$$\displaystyle
\mbox{\boldmath $f$}(x,\tau)=-\frac{\tau}{\epsilon}\tilde{f}(\tau)\chi_B(x)\mbox{\boldmath $a$},
\tag {1.9}
$$
where
$$\displaystyle
\tilde{f}(\tau)=\int_0^T e^{-\tau t}f(t)dt.
\tag {1.10}
$$

Before describing the main result we introduce two conditions
(A.I) and (A.II) listed below:

(A.I) $\exists C>0$\,\,\, $\displaystyle\lambda(x)\ge \sqrt{\frac{\epsilon}{\mu}}+C$ for all $x\in\partial D$;

(A.II)  $\exists C>0$\,$\exists C'>0$\,\,\, $\displaystyle C'\le\lambda(x)\le \sqrt{\frac{\epsilon}{\mu}}+C$ for all $x\in\partial D$.

Define $\text{dist}\,(D,B)=\inf_{x\in D,\,y\in B}\,\vert x-y\vert$.

\proclaim{\noindent Theorem 1.1.}  Let $\mbox{\boldmath $a$}_j$,
$j=1,2$ be two linearly independent unit vectors. Let
$\mbox{\boldmath $J$}_j(x,t)=f(t)\chi_B(x)\mbox{\boldmath $a$}_j$
and $f\in C^1[0,\,T]$ with $f(0)=0$ satisfy
$$\displaystyle
\exists \gamma\in\Bbb R\,\,\liminf_{\tau\longrightarrow\infty}\tau^{\gamma}\vert\tilde{f}(\tau)\vert>0.
\tag {1.11}
$$
Let $\mbox{\boldmath $f$}_j$, $j=1,2$ denote the $\mbox{\boldmath $f$}$ given by (1.9) and (1.10)  with $f$ above
and $\mbox{\boldmath $a$}=\mbox{\boldmath $a$}_j$.

Then, we have:

(i) For all $T\le2\sqrt{\mu\epsilon}\text{dist}\,(D,B)$ it holds that

$$\displaystyle
\lim_{\tau\longrightarrow\infty}e^{\tau T}\sum_{j=1}^2I_{\mbox{\boldmath $f$}_j}(\tau,T)=0;
$$

(ii) if $\lambda$ satisfies (A.I), then for all
$T>2\sqrt{\mu\epsilon}\text{dist}\,(D,B)$ it holds that

$$\displaystyle
\lim_{\tau\longrightarrow\infty}e^{\tau T}\sum_{j=1}^2I_{\mbox{\boldmath $f$}_j}(\tau,T)=\infty;
$$

(iii) if $\lambda$ satisfies (A.II), then for all
$T>2\sqrt{\mu\epsilon}\text{dist}\,(D,B)$ it holds that

$$\displaystyle
\lim_{\tau\longrightarrow\infty}e^{\tau T}\sum_{j=1}^2I_{\mbox{\boldmath $f$}_j}(\tau,T)=-\infty.
$$

Moreover, in case of both (ii) and (iii) we have, for all $T>2\sqrt{\mu\epsilon}\text{dist}\,(D,B)$
$$
\displaystyle
\lim_{\tau\longrightarrow\infty}\frac{1}{\tau}
\log\left\vert
\sum_{j=1}^2I_{\mbox{\boldmath $f$}_j}(\tau,T)\right\vert
=-2\sqrt{\mu\epsilon}\text{dist}\,(D,B).
\tag {1.12}
$$

\endproclaim

Theorem 1.1 says that the $T$ in the problem should be all $T$ satisfying
$$\displaystyle
T>2\sqrt{\mu\epsilon}\,\text{dist}\,(D,B).
$$
In particular, if we have a known upper bound $M$ such that $M>\text{dist}\,(D,B)$,
then we can choose $T=2\sqrt{\mu\epsilon}M$ and know the exact value
of $\text{dist}\,(D,B)$ via formula (1.12).

Moreover, we have another characterization of $\text{dist}\,(D,B)$.
Let $M$ be the same positive constant as above.  Assume that $\lambda$ satisfies
(A.I) or (A.II).  Given $f\in C^1[0,\,2\sqrt{\mu\epsilon}M]$ with $f(0)=0$ define
$$\begin{array}{ll}
\displaystyle
\tilde{f}_T(\tau)=\int_0^Te^{-\tau t}f(t)dt, &
\displaystyle
0<T<2\sqrt{\mu\epsilon}M.
\end{array}
$$
We assume that (1.11) is satisfied with $\tilde{f}=\tilde{f}_{2\sqrt{\mu\epsilon}M}$.
Generate $\mbox{\boldmath $E$}$ and $\mbox{\boldmath $H$}$ over time interval
$]0,\,2\sqrt{\mu\epsilon}M[$ using $\mbox{\boldmath $J$}(x,t)=f(t)\chi_B(x)\mbox{\boldmath $a$}_j$.
Measure $\mbox{\boldmath $E$}$ on $B$ over time interval $]0,\,2\sqrt{\mu\epsilon}M[$.
Compute also $\mbox{\boldmath $E$}_0$ on $B$ for $\mbox{\boldmath $J$}$ above over the same time interval.
Using those fields, for each $T\in]0,\,2\sqrt{\mu\epsilon}M[$
compute $\mbox{\boldmath $W$}_e$ and $\mbox{\boldmath $V$}_e$ given by
(1.7) and (1.8), respectively.  Denote by $\mbox{\boldmath $f$}_j^T$ the
$\mbox{\boldmath $f$}$ given by (1.6) and compute $I_{\mbox{\boldmath $f$}_j^T}(\tau,T)$
given by (1.5) for $\mbox{\boldmath $f$}=\mbox{\boldmath $f$}_j^T$.

Then, it is easy to see that (1.11) is satisfied also with $\tilde{f}=\tilde{f}_T$ for each
$T\in\,]0,\,2\sqrt{\mu\epsilon}M[$.
Therefore (i)-(iii) in Theorem 1.1 yield the formula
$$\displaystyle
]0,\,2\sqrt{\mu\epsilon}\,\text{dist}\,(D,B)]
=\left\{
T\in\,]0,\,2\sqrt{\mu\epsilon}M[\,\,\left\vert\,\,
\lim_{\tau\longrightarrow\infty}e^{\tau T}\sum_{j=1}^2 I_{\mbox{\boldmath $f$}_j^T}(\tau,T)=0\right\}.
\right.
\tag {1.13}
$$

Formula (1.13) has a similarity in the style as that of the {\it
original}l enclosure method applied to, for example, the Laplace
equation, see (1.3) in \cite{I1}. See also \cite{IW} for similar
statements to (ii) and (iii) in Theorem 1.1 for  scalar wave
equations in the whole space. Note that
$\text{dist}\,(D,B)=d_{\partial D}(p)-\eta$ where $p$ and $\eta$
are the center and radius of $B$, respectively and $d_{\partial
D}(p)=\inf_{x\in\partial D}\vert x-p\vert$. Thus knowing
$\text{dist}\,(D,B)$ is equivalent to knowing $d_{\partial D}(p)$.

We think that (i), (ii) and (iii) implicitly represent the {\it
finite propagation property} of the wave governed by the Maxwell
system.  It means that if $T\le
2\sqrt{\mu\epsilon}\text{dist}\,(D,B)$, then one can not get a
qualitative information about the state of the surface of the
obstacle (by using the enclosure method).  This is consistent with
that the propagation speed of the Maxwell system is given by
$1/\sqrt{\mu\epsilon}$.  However, note that in the proof presented
in this paper we do not make use of the finite propagation
property.

The results as stated in (ii) and (iii) together with formula (1.12)
can be considered as an extension of the corresponding results in \cite{IE2}.
More precisely, therein we considered an inverse obstacle scattering problem for the wave governed by the classical
wave equation $\partial_t^2u-\triangle u=0$ outside an obstacle which we denote by $D$ again.
The wave $u$ as a solution of the equation
satisfies $\partial u/\partial{\mbox{\boldmath $\nu$}}-\gamma\partial_t u-\beta u=0$
on $\partial D\times\,]0,\,T[$.  It is assumed that $\gamma$ and $\beta$
are essentially bounded functions on $\partial D$ and $\gamma\ge 0$.
Using an indicator function which can be computed from a wave generated by
a single set of initial data and a special solution of the modified Helmholtz equation, we found that $\gamma\equiv 1$ is
the special value as same as $\lambda\equiv\sqrt{\epsilon/\mu}$ of (A.I) and (A.II).
It means that the indicator function
therein changes the behaviour according to whether $\gamma$ is greater or less than $1$.

For the computation of indicator function (1.5) we need $\mbox{\boldmath $E$}_0$ on $B$ over time interval
$]0,\,T[$.  This can be done by solving explicitly and analytically the initial value problem (1.4).

Here we present another way which introduces an indicator function
using an approximation of (1.8).

Let $\mbox{\boldmath $V$}_e^0$ be the weak solution of
$$\begin{array}{ll}
\displaystyle
\frac{1}{\mu\epsilon}\nabla\times\nabla\times\mbox{\boldmath $V$}
+\tau^2\mbox{\boldmath $V$}+\mbox{\boldmath $f$}(x,\tau)
=\mbox{\boldmath $0$}
& \text{in}\,\Bbb R^3.
\end{array}
\tag {1.14}
$$

Define another indicator function
$$
\displaystyle
\tilde{I}_{\mbox{\boldmath $f$}}(\tau,T)
=\int_B\mbox{\boldmath $f$}(x,\tau)\cdot(\mbox{\boldmath $W$}_e-\mbox{\boldmath $V$}_e^0)dx.
\tag {1.15}
$$

The {\it theoretical advantage} of this indicator function compared with (1.5) is
that one has no need of computing $\mbox{\boldmath $V$}_e$ which requires the space time computation
of $\mbox{\boldmath $E$}_0$ and $\mbox{\boldmath $H$}_0$ on $B$ over time interval $]0,\,T[$.
Instead, in (1.15) one can compute $\mbox{\boldmath $V$}_e^0$ on $B$
in advance by solving only (1.14).
Our result on this indicator function is the following.

\proclaim{\noindent Theorem 1.2.}
All the statements of Theorem 1 are valid if $I_{\mbox{\boldmath $f$}_j}(\tau,T)$
is replaced with $\tilde{I}_{\mbox{\boldmath $f$}_j}(\tau,T)$.

\endproclaim

A brief outline of this paper is as follows. Theorems 1.1 and 1.2
are proven in Section 2 by using Lemmas 2.1-2.3. Lemma 2.1 gives
lower and upper estimates for $I_{\mbox{\boldmath $f$}}(\tau)$ as
$\tau\longrightarrow\infty$ with a single $\mbox{\boldmath $f$}$
and the proof is described in Section 3. The proof is based on a
{\it rough} asymptotic formula of $I_{\mbox{\boldmath $f$}}(\tau)$
as $\tau\longrightarrow\infty$ which is proved in Subsection 3.2.
Lemma 2.2 is concerned with an estimation of  the sum of the two
indicator functions corresponding to two input current sources in
the case when $\lambda$ satisfies (A.I)  or (A.II). It is proved
in Section 4. The statements (ii) and (iii) in Theorem 1.1 are
direct consequences of Lemmas 2.1 and 2.2. Lemma 2.3 describes a
simple estimate of the absolute value of $I_{\mbox{\boldmath
$f$}}(\tau)$ as $\tau\longrightarrow\infty$ which needs for the
proof of (i) and (1.12)  in Theorem 1.1.  The proof is given in
Subsection 4.5. The final section is devoted to conclusions and
some of problems to be solved in the future.

\section{Proof of Theorems 1.1 and 1.2}

One can obtain immediately the validity of the statements (ii) and (iii) in Theorem 1.1 from the following two lemmas.

\proclaim{\noindent Lemma 2.1.}
We have, as $\tau\longrightarrow\infty$
$$
\displaystyle
\tilde{J}_{e}(\tau)+O(\tau^{-1}e^{-\tau T})
\le I_{\mbox{\boldmath $f$}}(\tau,T)
\tag {2.1}
$$
and
$$
\displaystyle
I_{\mbox{\boldmath $f$}}(\tau,T)
\le
\tilde{J}_{e}(\tau)+
\frac{\tau}{\epsilon}\int_{\partial D}\lambda
\left\vert
\mbox{\boldmath $\nu$}\times(\mbox{\boldmath $V$}_e\times\mbox{\boldmath $\nu$})
+\frac{1}{\lambda}\mbox{\boldmath $V$}_m\times\mbox{\boldmath $\nu$}\right\vert^2dS
+O(\tau^{-1}e^{-\tau T}),
\tag {2.2}
$$
where
$$\displaystyle
\mbox{\boldmath $V$}_m=\mbox{\boldmath $V$}_m(x,\tau)
=\int_0^Te^{-\tau t}\mbox{\boldmath $H$}_0(x,\tau)dt
\tag {2.3}
$$
and
$$
\displaystyle
\tilde{J}_e(\tau)
=\frac{1}{\mu\epsilon}\int_{\partial D}(\mbox{\boldmath $\nu$}\times\mbox{\boldmath $V$}_e)
\cdot\nabla\times\mbox{\boldmath $V$}_e\,dS-
\frac{\tau}{\epsilon}\int_{\partial D}\frac{1}{\lambda}\vert\mbox{\boldmath $V$}_m\times\mbox{\boldmath $\nu$}\vert^2dS.
\tag {2.4}
$$

\endproclaim

Let $\mbox{\boldmath $V$}_{e,\,j}$ and $\mbox{\boldmath $V$}_{m,\,j}$ with $j=1,2$
denote $\mbox{\boldmath $V$}_e$ and $\mbox{\boldmath $V$}_m$ given by (1.8)
and (2.3) using $\mbox{\boldmath $E$}_0$ and $\mbox{\boldmath $H$}_0$ with $\mbox{\boldmath $J$}=\mbox{\boldmath $J$}_j$.

\proclaim{\noindent Lemma 2.2.}
Let $\tilde{J}_{e,\,j}$ with $j=1,2$ denote the $\tilde{J}_e$ given by (2.4) for $\mbox{\boldmath $V$}_{e}
=\mbox{\boldmath $V$}_{e,\,j}$ and $\mbox{\boldmath $V$}_m=\mbox{\boldmath $V$}_{m,\,j}$.

(i) If $\lambda$ satisfies (A.I), then there exist positive
numbers $\rho$, $C'$ and $\tau_0$ such that, for all
$\tau\ge\tau_0$
$$\displaystyle
\sum_{j=1}^2\tilde{J}_{e,\,j}(\tau)\ge C'\tau^{-\rho}e^{-2\tau\sqrt{\mu\epsilon}\text{dist}\,(D,B)}
+O(\tau^{-1/2}e^{-\tau T}).
$$

(ii) If $\lambda$ satisfies (A.II), then there exist positive
numbers $\rho$, $C'$ and $\tau_0$ such that, for all
$\tau\ge\tau_0$
$$\begin{array}{l}
\displaystyle
\sum_{j=1}^2
\left(\tilde{J}_{e,\,j}(\tau)
+\frac{\tau}{\epsilon}\int_{\partial D}\lambda
\left\vert
\mbox{\boldmath $\nu$}\times(\mbox{\boldmath $V$}_{e,\,j}\times\mbox{\boldmath $\nu$})
+\frac{1}{\lambda}\mbox{\boldmath $V$}_{m,\,j}\times\mbox{\boldmath $\nu$}\right\vert^2dS
\right)\\
\\
\displaystyle
\le -C'\tau^{-\rho }e^{-2\tau\sqrt{\mu\epsilon}\text{dist}\,(D,B)}
+O(\tau^{-1/2}e^{-\tau T}).
\end{array}
$$

\endproclaim

For the proof of (i) and formula (1.12) in Theorem 1.1 we need the following lemma.

\proclaim{\noindent Lemma 2.3.}
We have, as $\tau\longrightarrow\infty$
$$
\displaystyle
I_{\mbox{\boldmath $f$}}(\tau,T)
=O\left(\tau^{-1}e^{2\tau\sqrt{\mu\epsilon}\,\eta}(\tilde{f}(\tau))^2\int_{\partial D}v^2 dS\right)
+O(\tau^{-1/2}e^{-\tau T}).
\tag {2.5}
$$

\endproclaim

Since $\tilde{f}(\tau)=O(\tau^{-3/2})$, from (2.5) we have, as $\tau\longrightarrow\infty$
$$\displaystyle
e^{\tau T}I_{\mbox{\boldmath $f$}}(\tau,T)
=O(\tau^{-4}e^{\tau(T-2\sqrt{\mu\epsilon}\,\text{dist}\,(D,B))})+O(\tau^{-1/2}).
$$
This yields (i) of Theorem 1.1.  Furthermore a combination of this and Lemma 2.1 gives (1.12).
Note that $\rho$ in (i) and (ii) of Lemma 2.2 can be chosen as $\rho=2\gamma+5$.  See the end of the proof of  (i)
in Subsection 4.3.  Note that $\gamma$ in (1.11) has to be $\gamma\ge 3/2$ by the asymptotics
of $\tilde{f}(\tau)$ mentioned above.

Theorem 1.2 is a transplantation of Theorem 1.1 via the following simple estimate:
$$
\displaystyle
\Vert
\mbox{\boldmath $V$}_e-\mbox{\boldmath $V$}_e^0
\Vert_{L^2(\Bbb R^3)}
=O(\tau^{-1}e^{-\tau T}).
\tag {2.6}
$$
See (3.40) in Subsection 3.2.
Since it is easy to see that $\tilde{f}(\tau)=O(\tau^{-3/2})$ and hence $\Vert\mbox{\boldmath $f$}(\,\cdot\,,\tau)\Vert_{L^2(B)}
=O(\tau^{-1/2})$.
This together with (2.6) gives
$$
\displaystyle
\tilde{I}_{\mbox{\boldmath $f$}}(\tau,T)
=I_{\mbox{\boldmath $f$}}(\tau,T)+O(\tau^{-3/2}e^{-\tau T}).
$$
Now from this and Lemmas 2.1-2.3 we obtain Theorem 1.2.

Thus everything is reduced to giving the proof of Lemmas 2.1-2.3.

\section{Proof of Lemma 2.1}

\subsection{Preliminaries}

Define
$$\begin{array}{c}
\displaystyle
\mbox{\boldmath $W$}_m(x,\tau)
=\int_0^T e^{-\tau t}\mbox{\boldmath $H$}(x,t)dt.
\end{array}
$$
From (1.1) and (1.2) we see that $\mbox{\boldmath $W$}_e$ and $\mbox{\boldmath $W$}_m$
satisfy
$$
\left\{
\begin{array}{ll}
\displaystyle
\nabla\times\mbox{\boldmath $W$}_e+\tau\mu \mbox{\boldmath $W$}_m=-e^{-\tau T}\mu\mbox{\boldmath $H$}(x,T)
& \text{in}\,\Bbb R^3\setminus\overline D,\\
\\
\displaystyle
\nabla\times\mbox{\boldmath $W$}_m-\tau\epsilon \mbox{\boldmath $W$}_e
-\frac{\epsilon}{\tau}\mbox{\boldmath $f$}(x,\tau)
=e^{-\tau T}\epsilon \mbox{\boldmath $E$}(x,T)
& \text{in}\,\Bbb R^3\setminus\overline D
\end{array}
\right.
\tag {3.1}
$$
and
$$\begin{array}{ll}
\displaystyle
\mbox{\boldmath $\nu$}\times\mbox{\boldmath $W$}_m
-\lambda\,\mbox{\boldmath $\nu$}\times(\mbox{\boldmath $W$}_e\times\mbox{\boldmath $\nu$})=\mbox{\boldmath $0$}
& \text{on}\,\partial D.
\end{array}
\tag {3.2}
$$

Note that (3.2) is equivalent to
$$\begin{array}{ll}
\displaystyle
\mbox{\boldmath $\nu$}\times(\mbox{\boldmath $W$}_m\times\mbox{\boldmath $\nu$})
+\lambda\,\mbox{\boldmath $\nu$}\times\mbox{\boldmath $W$}_e=\mbox{\boldmath $0$}
& \text{on}\,\partial D.
\end{array}
\tag {3.3}
$$

Taking the rotation of equations on (3.1), we obtain
$$
\left\{
\begin{array}{ll}
\displaystyle
\frac{1}{\mu\epsilon}\nabla\times\nabla\times\mbox{\boldmath $W$}_e
+\tau^2\mbox{\boldmath $W$}_e+\mbox{\boldmath $f$}(x,\tau)
=e^{-\tau T}\mbox{\boldmath $F$}_e(x,\tau)
& \text{in}\,\Bbb R^3\setminus\overline D,
\\
\\
\displaystyle
\frac{1}{\mu\epsilon}\nabla\times\nabla\times\mbox{\boldmath $W$}_m
+\tau^2\mbox{\boldmath $W$}_m
-\frac{1}{\tau\mu}\nabla\times\mbox{\boldmath $f$}(x,\tau)
=e^{-\tau T}
\mbox{\boldmath $F$}_m(x,\tau)
& \text{in}\,\Bbb R^3\setminus\overline D,
\end{array}
\right.
\tag {3.4}
$$
where
$$
\left\{
\begin{array}{l}
\displaystyle
\mbox{\boldmath $F$}_e(x,\tau)=-\left(\tau\mbox{\boldmath $E$}(x,T)+\frac{1}{\epsilon}\nabla\times\mbox{\boldmath $H$}(x,T)\right),
\\
\\
\displaystyle
\mbox{\boldmath $F$}_m(x,\tau)=-\left(\tau\mbox{\boldmath $H$}(x,T)-\frac{1}{\mu}\nabla\times\mbox{\boldmath $E$}(x,T)\right).
\end{array}
\right.
\tag {3.5}
$$

Vector valued functions $\mbox{\boldmath $V$}_e$ and $\mbox{\boldmath $V$}_m$ satisfy
$$
\left\{
\begin{array}{ll}
\displaystyle
\nabla\times\mbox{\boldmath $V$}_m-\tau\epsilon\mbox{\boldmath $V$}_e-\frac{\epsilon}{\tau}
\mbox{\boldmath $f$}
=e^{-\tau T}\epsilon\mbox{\boldmath $E$}_0(x,T)
& \text{in}\,\Bbb R^3,
\\
\\
\displaystyle
\nabla\times\mbox{\boldmath $V$}_e+\tau\mu\mbox{\boldmath $V$}_m
=-e^{-\tau T}\mu\mbox{\boldmath $H$}_0(x,T)
& \text{in}\,\Bbb R^3.
\end{array}
\right.
\tag {3.6}
$$
Taking the rotation of equations on (3.6), we obtain
$$
\left\{
\begin{array}{ll}
\displaystyle
\frac{1}{\mu\epsilon}\nabla\times\nabla\times\mbox{\boldmath $V$}_e
+\tau^2\mbox{\boldmath $V$}_e+\mbox{\boldmath $f$}(x,\tau)
=e^{-\tau T}\mbox{\boldmath $F$}_e^0(x,\tau)
& \text{in}\,\Bbb R^3,
\\
\\
\displaystyle
\frac{1}{\mu\epsilon}\nabla\times\nabla\times\mbox{\boldmath $V$}_m
+\tau^2\mbox{\boldmath $V$}_m
-\frac{1}{\tau\mu}\nabla\times\mbox{\boldmath $f$}(x,\tau)
=e^{-\tau T}
\mbox{\boldmath $F$}_m^0(x,\tau)
& \text{in}\,\Bbb R^3,
\end{array}
\right.
\tag {3.7}
$$
where
$$
\left\{
\begin{array}{l}
\displaystyle
\mbox{\boldmath $F$}_e^0(x,\tau)=-\left(\tau\mbox{\boldmath $E$}_0(x,T)+\frac{1}{\epsilon}\nabla\times\mbox{\boldmath $H$}_0(x,T)\right),
\\
\\
\displaystyle
\mbox{\boldmath $F$}_m^0(x,\tau)=-\left(\tau\mbox{\boldmath $H$}_0(x,T)-\frac{1}{\mu}\nabla\times\mbox{\boldmath $E$}_0(x,T)\right).
\end{array}
\right.
\tag {3.8}
$$

Define
$$
\left\{
\begin{array}{c}
\displaystyle
\mbox{\boldmath $R$}_e=\mbox{\boldmath $W$}_e-\mbox{\boldmath $V$}_e,
\\
\\
\displaystyle
\mbox{\boldmath $R$}_m=\mbox{\boldmath $W$}_m-\mbox{\boldmath $V$}_m.
\end{array}
\right.
$$

From (3.1) and (3.6) we see that $\mbox{\boldmath $R$}_e$ and $\mbox{\boldmath $R$}_m$ satisfy
$$
\left\{
\begin{array}{ll}
\displaystyle
\nabla\times\mbox{\boldmath $R$}_m-\tau\epsilon\mbox{\boldmath $R$}_e
=e^{-\tau T}\epsilon\mbox{\boldmath $F$}
& \text{in}\,\Bbb R^3\setminus\overline D,\\
\\
\displaystyle
\nabla\times\mbox{\boldmath $R$}_e+\tau\mu\mbox{\boldmath $R$}_m
=-e^{-\tau T}\mu\mbox{\boldmath $G$}
& \text{in}\,\Bbb R^3\setminus\overline D,
\end{array}
\right.
\tag {3.9}
$$
where
$$
\left\{
\begin{array}{c}
\displaystyle
\mbox{\boldmath $F$}=\mbox{\boldmath $E$}(x,T)-\mbox{\boldmath $E$}_0(x,T),
\\
\\
\displaystyle
\mbox{\boldmath $G$}=\mbox{\boldmath $H$}(x,T)-\mbox{\boldmath $H$}_0(x,T).
\end{array}
\right.
\tag {3.10}
$$
Taking the differences of (3.4) from (3.7), we see that
$\mbox{\boldmath $R$}_{\star}$ with $\star=e,m$ satisfy
$$\begin{array}{ll}
\displaystyle
\frac{1}{\mu\epsilon}\nabla\times\nabla\times\mbox{\boldmath $R$}_{\star}+
\tau^2
\mbox{\boldmath $R$}_{\star}=e^{-\tau T}(\mbox{\boldmath $F$}_{\star}(x,\tau)-
\mbox{\boldmath $F$}_{\star}^0(x,\tau))
& \text{in}\,\Bbb R^3\setminus\overline D.
\end{array}
\tag {3.11}
$$
Define
$$
\displaystyle
E_{e}(\tau)
=\frac{1}{\epsilon\mu}
\int_{\Bbb R^3\setminus\overline D}\vert\nabla\times\mbox{\boldmath $R$}_{e}\vert^2dx
+\tau^2\int_{\Bbb R^3\setminus\overline D}\vert\mbox{\boldmath $R$}_{e}\vert^2dx.
\tag {3.12}
$$

\subsection{Rough asymptotic formula of the indicator function}

We start with having the following
asymptotic formula of the indicator function.

\proclaim{\noindent Proposition 3.1.}
It holds that, as $\tau\longrightarrow\infty$
$$
\displaystyle
\int_{\Bbb R^3\setminus\overline D}\mbox{\boldmath $f$}(x,\tau)\cdot\mbox{\boldmath $R$}_edx
=
\tilde{J}_e(\tau)+\tilde{E}_e(\tau)
+O(\tau^{-1}e^{-\tau T}),
\tag {3.13}
$$
where $\tilde{J}_e(\tau)$ is given by (2.4) and
$$\displaystyle
\tilde{E}_e(\tau)
=E_e(\tau)+\frac{\tau}{\epsilon}\int_{\partial D}\frac{1}{\lambda}\vert\mbox{\boldmath $R$}_m\times\mbox{\boldmath $\nu$}\vert^2\,dS.
\tag {3.14}
$$

\endproclaim

{\it\noindent Proof.}
The proof is divided into three steps.

{\noindent\it Step 1.}
First we show that
$$\begin{array}{c}
\displaystyle
\int_{\Bbb R^3\setminus\overline D}\mbox{\boldmath $f$}(x,\tau)\cdot\mbox{\boldmath $R$}_e\,dx
=
\tilde{J}_e(\tau)+\tilde{E}_e(\tau)+e^{-\tau T}(R_1(\tau)-R_2(\tau)),
\end{array}
\tag {3.15}
$$
where
$$\displaystyle
R_1(\tau)
=\frac{1}{\epsilon}
\int_{\partial D}
\frac{1}{\lambda}
(\mbox{\boldmath $G$}-\mbox{\boldmath $H$}_0)\cdot(\mbox{\boldmath $V$}_m\times\mbox{\boldmath $\nu$}+
\mbox{\boldmath $R$}_m\times\mbox{\boldmath $\nu$})\,dS
\tag {3.16}
$$
and
$$
\displaystyle
R_2(\tau)=\int_{\Bbb R^3\setminus\overline D}
\left\{
\mbox{\boldmath $F$}_e(x,\tau)\cdot\mbox{\boldmath $R$}_e
-
(\mbox{\boldmath $F$}_e(x,\tau)-\mbox{\boldmath $F$}_e^0(x,\tau))\cdot\mbox{\boldmath $V$}_e
\right\}\,dx.
\tag {3.17}
$$

This is proved as follows.
Integration by parts gives
$$\begin{array}{l}
\displaystyle
\,\,\,\,\,\,
\int_{\Bbb R^3\setminus\overline D}\left\{(\nabla\times\nabla\times\mbox{\boldmath $W$}_e)\cdot\mbox{\boldmath $V$}_e
-(\nabla\times\nabla\times\mbox{\boldmath $V$}_e)\cdot\mbox{\boldmath $W$}_e\right\}dx\\
\\
\displaystyle
=\int_{\partial D}
\left\{(\mbox{\boldmath $\nu$}\times(\nabla\times\mbox{\boldmath $V$}_e))\cdot\mbox{\boldmath $W$}_e
-(\mbox{\boldmath $\nu$}\times(\nabla\times\mbox{\boldmath $W$}_e))\cdot\mbox{\boldmath $V$}_e
\right\}dS.
\end{array}
$$
We have
$$\begin{array}{c}
\displaystyle
(\mbox{\boldmath $\nu$}\times(\nabla\times\mbox{\boldmath $V$}_e))\cdot\mbox{\boldmath $W$}_e
=(\mbox{\boldmath $W$}_e\times\mbox{\boldmath $\nu$})\cdot\nabla\times\mbox{\boldmath $V$}_e
=-(\mbox{\boldmath $\nu$}\times\mbox{\boldmath $W$}_e)\cdot\nabla\times\mbox{\boldmath $V$}_e
\end{array}
$$
and
$$\begin{array}{ll}
\displaystyle
(\mbox{\boldmath $\nu$}\times(\nabla\times\mbox{\boldmath $W$}_e))\cdot\mbox{\boldmath $V$}_e
& = (\nabla\times\mbox{\boldmath $W$}_e)\times\mbox{\boldmath $V$}_e\cdot\mbox{\boldmath $\nu$}
\\
\\
\displaystyle
& = (\mbox{\boldmath $V$}_e\times\mbox{\boldmath $\nu$})\cdot
(\nabla\times\mbox{\boldmath $W$}_e)\\
\\
\displaystyle
& = -(\mbox{\boldmath $\nu$}\times\mbox{\boldmath $V$}_e)\cdot
(\nabla\times\mbox{\boldmath $W$}_e).
\end{array}
$$
Thus
$$\begin{array}{l}
\displaystyle
\,\,\,\,\,\,
\int_{\Bbb R^3\setminus\overline D}\left\{(\nabla\times\nabla\times\mbox{\boldmath $W$}_e)\cdot\mbox{\boldmath $V$}_e
-(\nabla\times\nabla\times\mbox{\boldmath $V$}_e)\cdot\mbox{\boldmath $W$}_e\right\}dx\\
\\
\displaystyle
=\int_{\partial D}\left(-(\mbox{\boldmath $\nu$}\times\mbox{\boldmath $W$}_e)\cdot\nabla\times\mbox{\boldmath $V$}_e
+(\mbox{\boldmath $\nu$}\times\mbox{\boldmath $V$}_e)\cdot
(\nabla\times\mbox{\boldmath $W$}_e\right)dS.
\end{array}
$$
Substituting the first equations on (3.4) and (3.7) into this left-hand side, we obtain
$$
\begin{array}{l}
\displaystyle
\,\,\,\,\,\,
\frac{1}{\mu\epsilon}\int_{\partial D}\left(-(\mbox{\boldmath $\nu$}\times\mbox{\boldmath $W$}_e)\cdot\nabla\times\mbox{\boldmath $V$}_e
+(\mbox{\boldmath $\nu$}\times\mbox{\boldmath $V$}_e)\cdot
\nabla\times\mbox{\boldmath $W$}_e\right)dS
\\
\\
\displaystyle
=\int_{\Bbb R^3\setminus\overline D}\mbox{\boldmath $f$}(x,\tau)\cdot\mbox{\boldmath $R$}_e\,dx\\
\\
\displaystyle
\,\,\,
+e^{-\tau T}\int_{\Bbb R^3\setminus\overline D}
\left\{
(\mbox{\boldmath $F$}_e(x,\tau)-\mbox{\boldmath $F$}_e^0(x,\tau))\cdot\mbox{\boldmath $V$}_e
-\mbox{\boldmath $F$}_e^0(x,\tau)\cdot\mbox{\boldmath $R$}_e
\right\}
dx.
\end{array}
\tag {3.18}
$$
Using $\mbox{\boldmath $W$}_e=\mbox{\boldmath $V$}_e+\mbox{\boldmath $R$}_e$, one has
$$
\begin{array}{c}
\displaystyle
\,\,\,\,\,\,
\int_{\partial D}\left(-(\mbox{\boldmath $\nu$}\times\mbox{\boldmath $W$}_e)\cdot\nabla\times\mbox{\boldmath $V$}_e
+(\mbox{\boldmath $\nu$}\times\mbox{\boldmath $V$}_e)\cdot
\nabla\times\mbox{\boldmath $W$}_e\right)dS
\\
\\
\displaystyle
=\int_{\partial D}\left(-(\mbox{\boldmath $\nu$}\times\mbox{\boldmath $R$}_e)
\cdot\nabla\times\mbox{\boldmath $V$}_e
+(\mbox{\boldmath $\nu$}\times\mbox{\boldmath $V$}_e)\cdot\nabla\times\mbox{\boldmath $R$}_e\right)dS.
\end{array}
\tag {3.19}
$$
Note that from (3.3), we have
$$\begin{array}{ll}
\displaystyle
\mbox{\boldmath $\nu$}\times(\mbox{\boldmath $R$}_m\times\mbox{\boldmath $\nu$})
+\lambda\,\mbox{\boldmath $\nu$}\times\mbox{\boldmath $R$}_e
=-\mbox{\boldmath $\nu$}\times(\mbox{\boldmath $V$}_m\times\mbox{\boldmath $\nu$})
-\lambda\,\mbox{\boldmath $\nu$}\times\mbox{\boldmath $V$}_e
& \text{on}\,\partial D,
\end{array}
\tag {3.20}
$$
and thus
$$\begin{array}{ll}
\displaystyle
\mbox{\boldmath $\nu$}\times\mbox{\boldmath $V$}_e
=-\mbox{\boldmath $\nu$}\times\mbox{\boldmath $R$}_e
-\frac{1}{\lambda}\left\{\mbox{\boldmath $\nu$}\times(\mbox{\boldmath $V$}_m\times\mbox{\boldmath $\nu$})
+\mbox{\boldmath $\nu$}\times(\mbox{\boldmath $R$}_m\times\mbox{\boldmath $\nu$})\right\}
& \text{on}\,\partial D
\end{array}
$$
or equivalently,
$$\begin{array}{ll}
\displaystyle
\mbox{\boldmath $\nu$}\times\mbox{\boldmath $R$}_e
=-\mbox{\boldmath $\nu$}\times\mbox{\boldmath $V$}_e
-\frac{1}{\lambda}\left\{\mbox{\boldmath $\nu$}\times(\mbox{\boldmath $V$}_m\times\mbox{\boldmath $\nu$})
+\mbox{\boldmath $\nu$}\times(\mbox{\boldmath $R$}_m\times\mbox{\boldmath $\nu$})\right\}
& \text{on}\,\partial D.
\end{array}
$$
From these we obtain
$$\begin{array}{ll}
\displaystyle
(\mbox{\boldmath $\nu$}\times\mbox{\boldmath $V$}_e)\cdot
\nabla\times\mbox{\boldmath $R$}_e
=-\mbox{\boldmath $\nu$}\times\mbox{\boldmath $R$}_e\cdot
\nabla\times\mbox{\boldmath $R$}_e\\
\\
\displaystyle
-\frac{1}{\lambda}\left\{\mbox{\boldmath $\nu$}\times(\mbox{\boldmath $V$}_m\times\mbox{\boldmath $\nu$})
+\mbox{\boldmath $\nu$}\times(\mbox{\boldmath $R$}_m\times\mbox{\boldmath $\nu$})\right\}
\cdot
\nabla\times\mbox{\boldmath $R$}_e
& \text{on}\,\partial D
\end{array}
$$
and
$$\begin{array}{ll}
\displaystyle
(\mbox{\boldmath $\nu$}\times\mbox{\boldmath $R$}_e)\cdot
\nabla\times\mbox{\boldmath $V$}_e
=-\mbox{\boldmath $\nu$}\times\mbox{\boldmath $V$}_e\cdot
\nabla\times\mbox{\boldmath $V$}_e\\
\\
\displaystyle
-\frac{1}{\lambda}\left\{\mbox{\boldmath $\nu$}\times(\mbox{\boldmath $V$}_m\times\mbox{\boldmath $\nu$})
+\mbox{\boldmath $\nu$}\times(\mbox{\boldmath $R$}_m\times\mbox{\boldmath $\nu$})\right\}
\cdot
\nabla\times\mbox{\boldmath $V$}_e
& \text{on}\,\partial D.
\end{array}
$$
Therefore we obtain
$$
\begin{array}{l}
\displaystyle
\,\,\,\,\,\,
-(\mbox{\boldmath $\nu$}\times\mbox{\boldmath $R$}_e)\cdot
\nabla\times\mbox{\boldmath $V$}_e+(\mbox{\boldmath $\nu$}\times\mbox{\boldmath $V$}_e)\cdot
\nabla\times\mbox{\boldmath $R$}_e
\\
\\
\displaystyle
=\mbox{\boldmath $\nu$}\times\mbox{\boldmath $V$}_e\cdot
\nabla\times\mbox{\boldmath $V$}_e-\mbox{\boldmath $\nu$}\times\mbox{\boldmath $R$}_e\cdot
\nabla\times\mbox{\boldmath $R$}_e\\
\\
\displaystyle
\,\,\,
+\frac{1}{\lambda}\left\{\mbox{\boldmath $\nu$}\times(\mbox{\boldmath $V$}_m\times\mbox{\boldmath $\nu$})
+\mbox{\boldmath $\nu$}\times(\mbox{\boldmath $R$}_m\times\mbox{\boldmath $\nu$})\right\}\cdot
\nabla\times\mbox{\boldmath $V$}_e\\
\\
\displaystyle
\,\,\,
-\frac{1}{\lambda}\left\{\mbox{\boldmath $\nu$}\times(\mbox{\boldmath $V$}_m\times\mbox{\boldmath $\nu$})
+\mbox{\boldmath $\nu$}\times(\mbox{\boldmath $R$}_m\times\mbox{\boldmath $\nu$})\right\}
\cdot
\nabla\times\mbox{\boldmath $R$}_e.
\end{array}
\tag {3.21}
$$
Now from (3.18), (3.19) and (3.21) we obtain
$$
\begin{array}{l}
\displaystyle
\,\,\,\,\,\,
\int_{\Bbb R^3\setminus\overline D}\mbox{\boldmath $f$}(x,\tau)\cdot\mbox{\boldmath $R$}_e\,dx\\
\\
\displaystyle
=\frac{1}{\mu\epsilon}
\int_{\partial D}(\mbox{\boldmath $\nu$}\times\mbox{\boldmath $V$}_e\cdot
\nabla\times\mbox{\boldmath $V$}_e-\mbox{\boldmath $\nu$}\times\mbox{\boldmath $R$}_e\cdot
\nabla\times\mbox{\boldmath $R$}_e\,)\,dx\\
\\
\displaystyle\,\,\,
+\frac{1}{\mu\epsilon}
\int_{\partial D}
\frac{1}{\lambda}\left\{\mbox{\boldmath $\nu$}\times(\mbox{\boldmath $V$}_m\times\mbox{\boldmath $\nu$})
+\mbox{\boldmath $\nu$}\times(\mbox{\boldmath $R$}_m\times\mbox{\boldmath $\nu$})\right\}\cdot
\nabla\times\mbox{\boldmath $V$}_e\,dS\\
\\
\displaystyle
\,\,\,
-\frac{1}{\mu\epsilon}\int_{\partial D}\frac{1}{\lambda}\left\{\mbox{\boldmath $\nu$}\times(\mbox{\boldmath $V$}_m\times\mbox{\boldmath $\nu$})
+\mbox{\boldmath $\nu$}\times(\mbox{\boldmath $R$}_m\times\mbox{\boldmath $\nu$})\right\}
\cdot
\nabla\times\mbox{\boldmath $R$}_e\,dS
\\
\\
\displaystyle
\,\,\,
-e^{-\tau T}
\int_{\Bbb R^3\setminus\overline D}
\left\{
(\mbox{\boldmath $F$}_e(x,\tau)-\mbox{\boldmath $F$}_e^0(x,\tau))\cdot\mbox{\boldmath $V$}_e
-\mbox{\boldmath $F$}_e^0(x,\tau)\cdot\mbox{\boldmath $R$}_e\right\}
dx.
\end{array}
\tag {3.22}
$$
(3.11) with $\star=e$ and (3.12)  and integration by parts give
$$\begin{array}{l}
\displaystyle
\,\,\,\,\,\,
e^{-\tau T}\int_{\Bbb R^3\setminus\overline D}(\mbox{\boldmath $F$}_e(x,\tau)-\mbox{\boldmath $F$}_e^0(x,\tau))\cdot\mbox{\boldmath $R$}_e\,dx
\\
\\
\displaystyle
=-\frac{1}{\mu\epsilon}\int_{\partial D}\mbox{\boldmath $\nu$}\times(\nabla\times\mbox{\boldmath $R$}_e)
\cdot\mbox{\boldmath $R$}_e\,dS
+E_e(\tau)\\
\\
\displaystyle
=
\frac{1}{\mu\epsilon}\int_{\partial D}(\mbox{\boldmath $\nu$}\times\mbox{\boldmath $R$}_e)
\cdot\nabla\times\mbox{\boldmath $R$}_e\,dS+E_e(\tau),
\end{array}
$$
that is,
$$
\begin{array}{l}
\displaystyle
\,\,\,\,\,\,
-\frac{1}{\mu\epsilon}\int_{\partial D}(\mbox{\boldmath $\nu$}\times\mbox{\boldmath $R$}_e)
\cdot\nabla\times\mbox{\boldmath $R$}_e\,dS\\
\\
\displaystyle
=E_e(\tau)
-e^{-\tau T}\int_{\Bbb R^3\setminus\overline D}(\mbox{\boldmath $F$}_e(x,\tau)
-\mbox{\boldmath $F$}_e^0(x,\tau))\cdot\mbox{\boldmath $R$}_e\,dx.
\end{array}
\tag {3.23}
$$
Define
$$
\displaystyle
J_{e}(\tau)
=\frac{1}{\epsilon\mu}
\int_D\vert\nabla\times\mbox{\boldmath $V$}_{e}\vert^2dx
+\tau^2\int_{D}\vert\mbox{\boldmath $V$}_{e}\vert^2dx.
\tag {3.24}
$$
We have also from the first equation on (3.7) and (3.24)
$$
\begin{array}{c}
\displaystyle
\frac{1}{\mu\epsilon}\int_{\partial D}(\mbox{\boldmath $\nu$}\times\mbox{\boldmath $V$}_e)
\cdot\nabla\times\mbox{\boldmath $V$}_e\,dS
=J_e(\tau)-e^{-\tau T}\int_D\mbox{\boldmath $F$}_e^0(x,\tau)\cdot\mbox{\boldmath $V$}_e\,dx.
\end{array}
\tag {3.25}
$$
Now from (2.19), (2.20), (3.22) and (3.25), we obtain
$$
\begin{array}{l}
\displaystyle
\,\,\,\,\,\,
\int_{\Bbb R^3\setminus\overline D}\mbox{\boldmath $f$}(x,\tau)\cdot\mbox{\boldmath $R$}_e\,dx
\\
\\
\displaystyle
=J_e(\tau)-e^{-\tau T}
\int_D\mbox{\boldmath $F$}_e^0(x,\tau)\cdot\mbox{\boldmath $V$}_e\,dx+E_e(\tau)\\
\\
\displaystyle
\,\,\,
+\frac{1}{\mu\epsilon}\int_{\partial D}\frac{1}{\lambda}\left\{
\mbox{\boldmath $\nu$}\times(\mbox{\boldmath $R$}_m\times\mbox{\boldmath $\nu$})\cdot\nabla\times\mbox{\boldmath $V$}_e
-\mbox{\boldmath $\nu$}\times(\mbox{\boldmath $V$}_m\times\mbox{\boldmath $\nu$})\cdot
\nabla\times\mbox{\boldmath $R$}_e\right\}dS\\
\\
\displaystyle
\,\,\,
+\frac{1}{\mu\epsilon}\int_{\partial D}\frac{1}{\lambda}\left\{
-\mbox{\boldmath $\nu$}\times(\mbox{\boldmath $R$}_m\times\mbox{\boldmath $\nu$})
\cdot
\nabla\times\mbox{\boldmath $R$}_e
+(\mbox{\boldmath $\nu$}\times(\mbox{\boldmath $V$}_m\times\mbox{\boldmath $\nu$}))\cdot\nabla\times\mbox{\boldmath $V$}_e\right\}\,dS
\\
\\
\displaystyle
\,\,\,
-e^{-\tau T}\int_{\Bbb R^3\setminus\overline D}
\left\{
\mbox{\boldmath $F$}_e(x,\tau)\cdot\mbox{\boldmath $R$}_e
-
(\mbox{\boldmath $F$}_e(x,\tau)-\mbox{\boldmath $F$}_e^0(x,\tau))\cdot\mbox{\boldmath $V$}_e
\right\}\,dx.
\end{array}
\tag {3.26}
$$
Here we make an order of the integrals over $\partial D$ in (3.26).
It follows from the second equation on (3.6) that
$$
\begin{array}{l}
\displaystyle
\,\,\,\,\,\,
\mbox{\boldmath $\nu$}\times(\mbox{\boldmath $R$}_m\times\mbox{\boldmath $\nu$})\cdot\nabla\times\mbox{\boldmath $V$}_e
\\
\\
\displaystyle
=(\nabla\times\mbox{\boldmath $V$}_e)\times\mbox{\boldmath $\nu$}
\cdot(\mbox{\boldmath $R$}_m\times\mbox{\boldmath $\nu$})
\\
\\
\displaystyle
=-\tau\mu\mbox{\boldmath $V$}_m\times\mbox{\boldmath $\nu$}\cdot
\mbox{\boldmath $R$}_m\times\mbox{\boldmath $\nu$}
-e^{-\tau T}\mu\mbox{\boldmath $H$}_0\cdot(\mbox{\boldmath $R$}_m\times\mbox{\boldmath $\nu$})
\end{array}
\tag {3.27}
$$
and
$$
\begin{array}{l}
\displaystyle
\,\,\,\,\,\,
(\mbox{\boldmath $\nu$}\times(\mbox{\boldmath $V$}_m\times\mbox{\boldmath $\nu$}))\cdot\nabla\times\mbox{\boldmath $V$}_e
\\
\\
\displaystyle
=(\nabla\times\mbox{\boldmath $V$}_e)\times\mbox{\boldmath $\nu$}\cdot(\mbox{\boldmath $V$}_m\times\mbox{\boldmath $\nu$})
\\
\\
\displaystyle
=-\tau\mu\vert\mbox{\boldmath $V$}_m\times\mbox{\boldmath $\nu$}\vert^2
-e^{-\tau T}\mu\mbox{\boldmath $H$}_0\cdot(\mbox{\boldmath $V$}_m\times\mbox{\boldmath $\nu$}).
\end{array}
\tag {3.28}
$$

It follows from the second equation on (3.9) that
$$
\begin{array}{l}
\displaystyle
\,\,\,\,\,\,
\mbox{\boldmath $\nu$}\times(\mbox{\boldmath $V$}_m\times\mbox{\boldmath $\nu$})\cdot
\nabla\times\mbox{\boldmath $R$}_e
\\
\\
\displaystyle

=
(\nabla\times\mbox{\boldmath $R$}_e)\times\mbox{\boldmath $\nu$}
\cdot(\mbox{\boldmath $V$}_m\times\mbox{\boldmath $\nu$})
\\
\\
\displaystyle
=-\tau\mu\mbox{\boldmath $R$}_m\times\mbox{\boldmath $\nu$}\cdot\mbox{\boldmath $V$}_m\times\mbox{\boldmath $\nu$}
-e^{-\tau T}\mu\mbox{\boldmath $G$}\cdot(\mbox{\boldmath $V$}_m\times\mbox{\boldmath $\nu$})
\end{array}
\tag {3.29}
$$
and
$$
\begin{array}{l}
\displaystyle
\,\,\,\,\,\,
-(\mbox{\boldmath $\nu$}\times(\mbox{\boldmath $R$}_m\times\mbox{\boldmath $\nu$}))\cdot
\nabla\times\mbox{\boldmath $R$}_e
\\
\\
\displaystyle
=-(\nabla\times\mbox{\boldmath $R$}_e)\times\mbox{\boldmath $\nu$}
\cdot(\mbox{\boldmath $R$}_m\times\mbox{\boldmath $\nu$})\\
\\
\displaystyle
=\tau\mu\vert\mbox{\boldmath $R$}_m\times\mbox{\boldmath $\nu$}\vert^2+
e^{-\tau T}\mu\mbox{\boldmath $G$}\cdot(\mbox{\boldmath $R$}_m\times\mbox{\boldmath $\nu$}).
\end{array}
\tag {3.30}
$$
A combination of (3.27) and (3.29) gives
$$
\begin{array}{l}
\displaystyle
\,\,\,\,\,\,
\mbox{\boldmath $\nu$}\times(\mbox{\boldmath $R$}_m\times\mbox{\boldmath $\nu$})\cdot\nabla\times\mbox{\boldmath $V$}_e
-\mbox{\boldmath $\nu$}\times(\mbox{\boldmath $V$}_m\times\mbox{\boldmath $\nu$})\cdot
\nabla\times\mbox{\boldmath $R$}_e\\
\\
\displaystyle
=e^{-\tau T}\mu
\left\{\mbox{\boldmath $G$}\cdot(\mbox{\boldmath $V$}_m\times\mbox{\boldmath $\nu$})
-\mbox{\boldmath $H$}_0\cdot(\mbox{\boldmath $R$}_m\times\mbox{\boldmath $\nu$})
\right\}.
\end{array}
\tag {3.31}
$$
A combination of (3.28) and (3.30) gives
$$
\begin{array}{l}
\displaystyle
\,\,\,\,\,\,
-(\mbox{\boldmath $\nu$}\times(\mbox{\boldmath $R$}_m\times\mbox{\boldmath $\nu$}))
\cdot
\nabla\times\mbox{\boldmath $R$}_e
+(\mbox{\boldmath $\nu$}\times(\mbox{\boldmath $V$}_m\times\mbox{\boldmath $\nu$}))\cdot\nabla\times\mbox{\boldmath $V$}_e\\
\\
\displaystyle
=-\tau\mu\vert\mbox{\boldmath $V$}_m\times\mbox{\boldmath $\nu$}\vert^2+\tau\mu\vert\mbox{\boldmath $R$}_m\times\mbox{\boldmath $\nu$}\vert^2
+e^{-\tau T}\mu\left\{\mbox{\boldmath $G$}\cdot(\mbox{\boldmath $R$}_m\times\mbox{\boldmath $\nu$})
-\mbox{\boldmath $H$}_0\cdot(\mbox{\boldmath $V$}_m\times\mbox{\boldmath $\nu$})\right\}.
\end{array}
\tag {3.32}
$$
It follows from (2.4) and (3.25) that
$$\displaystyle
\tilde{J}_e(\tau)
=
J_e(\tau)-\frac{\tau}{\epsilon}\int_{\partial D}\frac{1}{\lambda}\vert\mbox{\boldmath $V$}_m\times\mbox{\boldmath $\nu$}\vert^2dS
-e^{-\tau T}
\int_D\mbox{\boldmath $F$}_e^0(x,\tau)\cdot\mbox{\boldmath $V$}_e\,dx.
\tag {3.33}
$$
Now substituting (3.31) , (3.32)  and (3.33) into (3.26), we obtain (3.15).

{\it\noindent Step 2.}  We have, as $\longrightarrow\infty$

$$\displaystyle
\Vert\mbox{\boldmath $V$}_e\Vert_{L^2(\Bbb R^3)}=O(\tau^{-5/2});
\tag {3.34}
$$
for $\star=e,m$
$$\displaystyle
\Vert\mbox{\boldmath $V$}_{\star}\times\mbox{\boldmath $\nu$}\Vert_{L^2(\partial D)}=O(\tau^{-3/2}).
\tag {3.35}
$$

This is proved as follows.

Recall $\mbox{\boldmath $V$}_e^0$ is the weak solution of (1.14).
Since $\Vert \mbox{\boldmath $f$}(\,\cdot\,\tau)\Vert_{L^2(B)}=O(\tau^{-1/2})$,
it is easy to see that we have
$$
\Vert\mbox{\boldmath $V$}_e^0\Vert_{L^2(\Bbb R^3)}
=O(\tau^{-5/2}).
\tag {3.36}
$$
and
$$\displaystyle
\Vert
\nabla\times\mbox{\boldmath $V$}_e^0
\Vert_{L^2(\Bbb R^3)}
=O(\tau^{-3/2}).
\tag {3.37}
$$

Next define
$$\displaystyle
\mbox{\boldmath $Z$}=\mbox{\boldmath $V$}_e-\mbox{\boldmath $V$}_e^0.
$$
Then from the first equation on (3.7) and (1.14) we have
$$
\begin{array}{ll}
\displaystyle
\frac{1}{\mu\epsilon}\nabla\times\nabla\times\mbox{\boldmath $Z$}
+\tau^2\mbox{\boldmath $Z$}=e^{-\tau T}\mbox{\boldmath $F$}_e^0
& \text{in}\,\Bbb R^n.
\end{array}
\tag {3.38}
$$
This gives
$$\displaystyle
\frac{1}{\mu\epsilon}\int_{\Bbb R^n}\vert\nabla\times\mbox{\boldmath $Z$}\vert^2dx
+\tau^2\int_{\Bbb R^3}\vert\mbox{\boldmath $Z$}\vert^2dx
=e^{-\tau T}\int_{\Bbb R^3}\mbox{\boldmath $F$}_e^0\cdot\mbox{\boldmath $Z$}dx
$$
and hence
$$\displaystyle
\frac{1}{\mu\epsilon}\int_{\Bbb R^3}\vert\nabla\times\mbox{\boldmath $Z$}\vert^2dx
+\tau^2\int_{\Bbb R^3}\vert\mbox{\boldmath $Z$}\vert^2dx
\le
\frac{e^{-2\tau T}}{\tau^2}\int_{\Bbb R^3}\vert\mbox{\boldmath $F$}_e^0\vert^2dx.
$$
Then applying the first equation on (3.8) to this right-hand side, we conclude that
$$
\displaystyle\Vert\mbox{\boldmath $Z$}\Vert_{L^2(\Bbb R^3)}=O(\tau^{-1}e^{-\tau T}),
\tag {3.39}
$$
and
$$
\displaystyle\Vert\nabla\times\mbox{\boldmath $Z$}\Vert_{L^2(\Bbb R^3)}=O(e^{-\tau T}).
\tag {3.40}
$$
Now a combination of (3.36) and (3.39) gives (3.34).
Note also that a combination of  (3.37) and (3.40) gives
$$
\displaystyle\Vert\nabla\times\mbox{\boldmath $V$}_e\Vert_{L^2(\Bbb R^3)}=O(\tau^{-3/2}).
\tag {3.41}
$$
Then, the trace theorem \cite{DL3} yields (3.35) with $\star=e$.
Moreover, using equations (3.6) together with (3.34) and (3.41), we obtain
$$
\displaystyle\Vert\mbox{\boldmath $V$}_m\Vert_{L^2(\Bbb R^3)}=O(\tau^{-5/2}),
$$
and
$$
\displaystyle\Vert\nabla\times\mbox{\boldmath $V$}_m\Vert_{L^2(\Bbb R^3)}=O(\tau^{-3/2}).
$$
Then, the trace theorem \cite{DL3} yields (3.35) with $\star=m$.

{\it\noindent Step 3.} We have, as $\tau\longrightarrow\infty$
$$
\displaystyle
\Vert\mbox{\boldmath $R$}_e\Vert_{L^2(\Bbb R^3\setminus\overline D)}=O(\tau^{-2});
\tag {3.42}
$$
$$
\displaystyle
\Vert\mbox{\boldmath $R$}_m\times\mbox{\boldmath $\nu$}\Vert_{L^2(\partial D)}
=O(\tau^{-3/2}).
\tag {3.43}
$$

This is proved as follows.

From (3.9) we obtain
$$\displaystyle
\begin{array}{l}
\,\,\,\,\,\,
\displaystyle
\int_{\Bbb R^3\setminus\overline D}\nabla\cdot(\mbox{\boldmath $R$}_e\times\mbox{\boldmath $R$}_m)dx
+
\tau\int_{\Bbb R^3\setminus\overline D}(\epsilon\vert\mbox{\boldmath $R$}_e\vert^2+\mu\vert\mbox{\boldmath $R$}_m\vert^2)dx\\
\\
\displaystyle
=-e^{-\tau T}\int_{\Bbb R^3\setminus\overline D}(\epsilon\mbox{\boldmath $F$}\cdot\mbox{\boldmath $R$}_e
+\mu\mbox{\boldmath $G$}\cdot\mbox{\boldmath $R$}_m)dx.
\end{array}
\tag {3.44}
$$
Here we note that
$$\displaystyle
\int_{\Bbb R^3\setminus\overline D}\nabla\cdot(\mbox{\boldmath $R$}_e\times\mbox{\boldmath $R$}_m)dx
=-\int_{\partial D}\mbox{\boldmath $\nu$}\cdot\mbox{\boldmath $R$}_e\times\mbox{\boldmath $R$}_mdS
$$
and from (3.20) we have
$$\displaystyle
\begin{array}{l}
\displaystyle
\,\,\,\,\,\,
\mbox{\boldmath $\nu$}\cdot\mbox{\boldmath $R$}_e\times\mbox{\boldmath $R$}_m
\\
\\
\displaystyle
=\mbox{\boldmath $R$}_m\cdot(\mbox{\boldmath $\nu$}\times\mbox{\boldmath $R$}_e)\\
\\
\displaystyle
=\mbox{\boldmath $\nu$}\times(\mbox{\boldmath $R$}_m\times\mbox{\boldmath $\nu$})
\cdot
(\mbox{\boldmath $\nu$}\times\mbox{\boldmath $R$}_e)\\
\\
\displaystyle
=-\lambda\mbox{\boldmath $\nu$}\times\mbox{\boldmath $R$}_e\cdot\mbox{\boldmath $\nu$}\times\mbox{\boldmath $R$}_e
-\mbox{\boldmath $\nu$}\times(\mbox{\boldmath $V$}_m\times\mbox{\boldmath $\nu$})
\cdot
\mbox{\boldmath $\nu$}\times\mbox{\boldmath $R$}_e
-\lambda\mbox{\boldmath $\nu$}\times\mbox{\boldmath $V$}_e\cdot
\mbox{\boldmath $\nu$}\times\mbox{\boldmath $R$}_e\\
\\
\displaystyle
=-\lambda\vert\mbox{\boldmath $\nu$}\times\mbox{\boldmath $R$}_e\vert^2
-(\mbox{\boldmath $\nu$}\times(\mbox{\boldmath $V$}_m\times\mbox{\boldmath $\nu$})+\lambda\mbox{\boldmath $\nu$}\times\mbox{\boldmath $V$}_e)
\cdot\mbox{\boldmath $\nu$}\times\mbox{\boldmath $R$}_e\\
\\
\displaystyle
=-\lambda\left\vert\mbox{\boldmath $\nu$}\times\mbox{\boldmath $R$}_e+
\frac{1}{2}
\left(\mbox{\boldmath $\nu$}\times\mbox{\boldmath $V$}_e+\frac{1}{\lambda}
\mbox{\boldmath $\nu$}\times(\mbox{\boldmath $V$}_m\times\mbox{\boldmath $\nu$})
\right)
\right\vert^2\\
\\
\displaystyle
\,\,\,
+\frac{\lambda}{4}
\left\vert\mbox{\boldmath $\nu$}\times\mbox{\boldmath $V$}_e+\frac{1}{\lambda}
\mbox{\boldmath $\nu$}\times(\mbox{\boldmath $V$}_m\times\mbox{\boldmath $\nu$})
\right\vert^2.
\end{array}
\tag {3.45}
$$
Therefore (3.44) becomes
$$\begin{array}{l}
\displaystyle
\,\,\,\,\,\,
\int_{\partial D}\lambda\left\vert\mbox{\boldmath $\nu$}\times\mbox{\boldmath $R$}_e+
\frac{1}{2}
\left(\mbox{\boldmath $\nu$}\times\mbox{\boldmath $V$}_e+\frac{1}{\lambda}
\mbox{\boldmath $\nu$}\times(\mbox{\boldmath $V$}_m\times\mbox{\boldmath $\nu$})
\right)
\right\vert^2 dS
+\tau\int_{\Bbb R^3\setminus\overline D}(\epsilon\vert\mbox{\boldmath $R$}_e\vert^2+\mu\vert\mbox{\boldmath $R$}_m\vert^2)dx\\
\\
\displaystyle
\,\,\,
+e^{-\tau T}\int_{\Bbb R^3\setminus\overline D}(\epsilon\mbox{\boldmath $F$}\cdot\mbox{\boldmath $R$}_e
+\mu\mbox{\boldmath $G$}\cdot\mbox{\boldmath $R$}_m)dx
\\
\\
\displaystyle
=\frac{1}{4}\int_{\partial D}\lambda
\left\vert\mbox{\boldmath $\nu$}\times\mbox{\boldmath $V$}_e+\frac{1}{\lambda}
\mbox{\boldmath $\nu$}\times(\mbox{\boldmath $V$}_m\times\mbox{\boldmath $\nu$})
\right\vert^2 dS
\end{array}
$$
and hence
$$
\begin{array}{l}
\displaystyle
\,\,\,\,\,\,
\int_{\partial D}\lambda\left\vert\mbox{\boldmath $\nu$}\times\mbox{\boldmath $R$}_e+
\frac{1}{2}
\left(\mbox{\boldmath $\nu$}\times\mbox{\boldmath $V$}_e+\frac{1}{\lambda}
\mbox{\boldmath $\nu$}\times(\mbox{\boldmath $V$}_m\times\mbox{\boldmath $\nu$})
\right)
\right\vert^2 dS\\
\\
\displaystyle
\,\,\,
+
\tau
\int_{\Bbb R^3\setminus\overline D}
\left(
\epsilon
\left\vert\mbox{\boldmath $R$}_e+\frac{e^{-\tau T}}{2\tau}\mbox{\boldmath $F$}\right\vert^2
+\mu
\left\vert\mbox{\boldmath $R$}_m+\frac{e^{-\tau T}}{2\tau}\mbox{\boldmath $G$}\right\vert^2
\right)
dx\\
\\
\displaystyle
=\frac{1}{4}\int_{\partial D}\lambda
\left\vert\mbox{\boldmath $\nu$}\times\mbox{\boldmath $V$}_e+\frac{1}{\lambda}
\mbox{\boldmath $\nu$}\times(\mbox{\boldmath $V$}_m\times\mbox{\boldmath $\nu$})
\right\vert^2 dS
+\frac{e^{-2\tau T}}{4\tau}
\int_{\Bbb R^3\setminus\overline D}
\left(\epsilon\vert\mbox{\boldmath $F$}\vert^2+
\mu\vert\mbox{\boldmath $G$}\vert^2
\right)dx.
\end{array}
$$
This yields
$$\begin{array}{l}
\displaystyle
\,\,\,\,\,\,
\frac{1}{2}
\int_{\partial D}\lambda\vert\mbox{\boldmath $\nu$}\times\mbox{\boldmath $R$}_e\vert^2 dS
+
\frac{1}{2}\tau
\int_{\Bbb R^3\setminus\overline D}
(\epsilon\vert\mbox{\boldmath $R$}_e\vert^2+\mu\vert\mbox{\boldmath $R$}_m\vert^2)dx
\\
\\
\displaystyle
\le
\frac{1}{2}
\int_{\partial D}
\lambda\left\vert
\mbox{\boldmath $\nu$}\times\mbox{\boldmath $V$}_e+
\frac{1}{\lambda}\mbox{\boldmath $\nu$}\times(\mbox{\boldmath $V$}_m\times\mbox{\boldmath $\nu$})
\right\vert^2
dS+
\frac{e^{-2\tau T}}{2\tau}
\int_{\Bbb R^3\setminus\overline D}(\epsilon\vert\mbox{\boldmath $F$}\vert^2+\mu\vert\mbox{\boldmath $G$}\vert^2)dx,
\end{array}
$$
and hence
$$
\begin{array}{l}
\displaystyle
\,\,\,\,\,\,
\lambda_0\int_{\partial D}\vert\mbox{\boldmath $\nu$}\times\mbox{\boldmath $R$}_e\vert^2 dS
+
\tau
\int_{\Bbb R^3\setminus\overline D}
(\epsilon\vert\mbox{\boldmath $R$}_e\vert^2+\mu\vert\mbox{\boldmath $R$}_m\vert^2)dx
\\
\\
\displaystyle
\le
\int_{\partial D}
\lambda\left\vert
\mbox{\boldmath $\nu$}\times\mbox{\boldmath $V$}_e+
\frac{1}{\lambda}
\mbox{\boldmath $\nu$}\times(\mbox{\boldmath $V$}_m\times\mbox{\boldmath $\nu$})
\right\vert^2
dS+
\frac{e^{-2\tau T}}{\tau}
\int_{\Bbb R^3\setminus\overline D}(\epsilon\vert\mbox{\boldmath $F$}\vert^2+\mu\vert\mbox{\boldmath $G$}\vert^2)dx,
\end{array}
\tag {3.46}
$$
where $\lambda_0=\min_{x\in\partial D}\lambda(x)$.
Now applying (3.35) and (3.10) to the right-hand side on (3.46)
we obtain (3.42).

The proof of estimate (3.43) is the following.
First from (3.2) we have
$$\begin{array}{c}
\displaystyle
\mbox{\boldmath $\nu$}\times(\mbox{\boldmath $R$}_e\times\mbox{\boldmath $\nu$})
=\frac{1}{\lambda}
\mbox{\boldmath $\nu$}\times\mbox{\boldmath $R$}_m
-\mbox{\boldmath $V$}_{em},
\end{array}
$$
where
$$\displaystyle
\mbox{\boldmath $V$}_{em}
=
\mbox{\boldmath $\nu$}\times(\mbox{\boldmath $V$}_e\times\mbox{\boldmath $\nu$})
+\frac{1}{\lambda}\mbox{\boldmath $V$}_m\times\mbox{\boldmath $\nu$}.
\tag {3.47}
$$
Thus we have another expression for (3.45):
$$\begin{array}{l}
\displaystyle
\,\,\,\,\,\,
\mbox{\boldmath $\nu$}\cdot\mbox{\boldmath $R$}_e\times\mbox{\boldmath $R$}_m
\\
\\
\displaystyle
=-\mbox{\boldmath $R$}_e\cdot(\mbox{\boldmath $\nu$}\times\mbox{\boldmath $R$}_m)\\
\\
\displaystyle
=-\mbox{\boldmath $\nu$}\times(\mbox{\boldmath $R$}_e\times\mbox{\boldmath $\nu$})
\cdot
(\mbox{\boldmath $\nu$}\times\mbox{\boldmath $R$}_m)\\
\\
\displaystyle
=-\frac{1}{\lambda}\vert\mbox{\boldmath $\nu$}\times\mbox{\boldmath $R$}_m\vert^2
+\mbox{\boldmath $V$}_{em}
\cdot\mbox{\boldmath $\nu$}\times\mbox{\boldmath $R$}_m\\
\\
\displaystyle
=-\frac{1}{\lambda}\left\vert\mbox{\boldmath $\nu$}\times\mbox{\boldmath $R$}_m-
\frac{\lambda}{2}\mbox{\boldmath $V$}_{em}
\right\vert^2
+\frac{\lambda}{4}
\left\vert
\mbox{\boldmath $V$}_{em}\right\vert^2.
\end{array}
$$
Therefore (3.44) becomes
$$\begin{array}{l}
\displaystyle
\,\,\,\,\,\,
\int_{\partial D}\frac{1}{\lambda}\left\vert\mbox{\boldmath $\nu$}\times\mbox{\boldmath $R$}_m-
\frac{\lambda}{2}
\mbox{\boldmath $V$}_{em}
\right\vert^2 dS
+\tau\int_{\Bbb R^3\setminus\overline D}(\epsilon\vert\mbox{\boldmath $R$}_e\vert^2+\mu\vert\mbox{\boldmath $R$}_m\vert^2)dx\\
\\
\displaystyle
\,\,\,
+e^{-\tau T}\int_{\Bbb R^3\setminus\overline D}(\epsilon\mbox{\boldmath $F$}\cdot\mbox{\boldmath $R$}_e
+\mu\mbox{\boldmath $G$}\cdot\mbox{\boldmath $R$}_m)dx
\\
\\
\displaystyle
=\frac{1}{4}\int_{\partial D}\lambda
\left\vert
\mbox{\boldmath $V$}_{em}\right\vert^2 dS
\end{array}
$$
and hence
$$
\begin{array}{l}
\displaystyle
\,\,\,\,\,\,
\int_{\partial D}\frac{1}{\lambda}\left\vert\mbox{\boldmath $\nu$}\times\mbox{\boldmath $R$}_m-
\frac{\lambda}{2}
\mbox{\boldmath $V$}_{em}
\right\vert^2 dS
\\
\\
\displaystyle
\,\,\,
+
\tau
\int_{\Bbb R^3\setminus\overline D}
\left(
\epsilon
\left\vert\mbox{\boldmath $R$}_e+\frac{e^{-\tau T}}{2\tau}\mbox{\boldmath $F$}\right\vert^2
+\mu
\left\vert\mbox{\boldmath $R$}_m+\frac{e^{-\tau T}}{2\tau}\mbox{\boldmath $G$}\right\vert^2
\right)
dx\\
\\
\displaystyle
=\frac{1}{4}\int_{\partial D}\lambda
\left\vert
\mbox{\boldmath $V$}_{em}\right\vert^2 dS
+\frac{e^{-2\tau T}}{4\tau}
\int_{\Bbb R^3\setminus\overline D}
\left(\epsilon\vert\mbox{\boldmath $F$}\vert^2+
\mu\vert\mbox{\boldmath $G$}\vert^2
\right)dx.
\end{array}
$$
This yields
$$
\begin{array}{l}
\displaystyle
\,\,\,\,\,\,
\frac{1}{\lambda_1}\int_{\partial D}\vert\mbox{\boldmath $\nu$}\times\mbox{\boldmath $R$}_m\vert^2 dS
+
\tau
\int_{\Bbb R^3\setminus\overline D}
(\epsilon\vert\mbox{\boldmath $R$}_e\vert^2+\mu\vert\mbox{\boldmath $R$}_m\vert^2)dx
\\
\\
\displaystyle
\le
\int_{\partial D}\lambda
\left\vert
\mbox{\boldmath $V$}_{em}\right\vert^2 dS
+
\frac{e^{-2\tau T}}{\tau}
\int_{\Bbb R^3\setminus\overline D}(\epsilon\vert\mbox{\boldmath $F$}\vert^2+\mu\vert\mbox{\boldmath $G$}\vert^2)dx,
\end{array}
\tag {3.48}
$$
where $\lambda_1=\max_{x\in\partial D}\lambda(x)$.
Now applying (3.35) and (3.10) to the right-hand side on (3.48),
we obtain (3.43).

Finally applying (3.35) and (3.43) to the right-hand side on  (3.16),
we obtain
$$
\displaystyle
R_1(\tau)
=O(\tau^{-3/2}).
$$
Applying (3.34) and (3.41) to the right-hand side on (3.17) and noting (3.5) and (3.8),
we obtain
$$
\displaystyle
R_2(\tau)
=O(\tau^{-1})+O(\tau^{-3/2})
=O(\tau^{-1}).
$$

\noindent
$\Box$

{\bf\noindent Remark 3.1.} Needless to say, estimates (3.35),
(3.42) and (3.43) are quite rough.  In fact, using the explicit
expression of $\mbox{\boldmath $V$}_e^0$ outside $B$, we see that
$\mbox{\boldmath $V$}_e^0$ together with its derivatives is
exponentially decaying as $\tau\longrightarrow\infty$ uniformly
for all $x\in\overline D$. This yields better estimates than
(3.35), (3.42) and (3.43). However, at this stage we do not need
such a detailed information about $\mbox{\boldmath $V$}_e^0$. This
is a reason why we call (3.13) the {\it rough} asymptotic formula.
Note that (3.39) is nothing but (2.6).

\subsection{Finishing the proof of Lemma 2.1}

In this subsection we derive the following estimate
for $\tilde{E}_e(\tau)$ given by (3.14):
$$
\begin{array}{c}
\displaystyle
\tilde{E}_e(\tau)
\le
\frac{\tau}{\epsilon}
\int_{\partial D}\lambda
\vert
\mbox{\boldmath $V$}_{em}\vert^2dS
+O(e^{-2\tau T}),
\end{array}
\tag {3.49}
$$
where $\mbox{\boldmath $V$}_{em}$ is given by (3.47).
Note that from this together with (3.13) we obtain (2.1) and (2.2).

It follows from (3.12) and (3.23) that
$$
\begin{array}{l}
\displaystyle
\,\,\,\,\,\,
\frac{1}{\mu\epsilon}\int_{\Bbb R^3\setminus\overline D}\vert\nabla\times\mbox{\boldmath $R$}_e\vert^2dx
+\tau^2\int_{\Bbb R^3\setminus\overline D}\left\vert\mbox{\boldmath $R$}_e-
\frac{e^{-\tau T}}{2\tau^2}(\mbox{\boldmath $F$}_e-\mbox{\boldmath $F$}_e^0)\right\vert^2dx\\
\\
\displaystyle
\,\,\,
+\frac{1}{\mu\epsilon}\int_{\partial D}(\mbox{\boldmath $\nu$}\times\mbox{\boldmath $R$}_e)
\cdot\nabla\times\mbox{\boldmath $R$}_e\,dS
\\
\\
\displaystyle
=\frac{e^{-2\tau T}}{4\tau^2}\int_{\Bbb R^3\setminus\overline D}\vert\mbox{\boldmath $F$}_e-\mbox{\boldmath $F$}_0\vert^2dx.
\end{array}
\tag {3.50}
$$
Since
$$\displaystyle
\left\vert
\mbox{\boldmath $R$}_e-\frac{e^{-\tau T}}{2\tau^2}(\mbox{\boldmath $F$}_e-\mbox{\boldmath $F$}_e^0)\right\vert^2
\ge\frac{1}{2}\vert\mbox{\boldmath $R$}_e\vert^2
-\frac{e^{-2\tau T}}{4\tau^4}\vert\mbox{\boldmath $F$}_e-\mbox{\boldmath $F$}_e^0\vert^2,
$$
It follows from (3.50) that
$$
\displaystyle
\frac{1}{2}E_e(\tau)+\frac{1}{\mu\epsilon}\int_{\partial D}(\mbox{\boldmath $\nu$}\times\mbox{\boldmath $R$}_e)
\cdot\nabla\times\mbox{\boldmath $R$}_e\,dS
\le
\frac{e^{-2\tau T}}{2\tau^2}\int_{\Bbb R^3\setminus\overline D}\vert\mbox{\boldmath $F$}_e-\mbox{\boldmath $F$}_0\vert^2dx.
\tag {3.51}
$$
From (3.3) we have
$$\begin{array}{l}
\displaystyle
\,\,\,\,\,\,
\mbox{\boldmath $\nu$}\times\mbox{\boldmath $R$}_e\\
\\
\displaystyle
=\mbox{\boldmath $\nu$}\times(\mbox{\boldmath $W$}_e-\mbox{\boldmath $V$}_e)\\
\\
\displaystyle
=-\frac{1}{\lambda}\mbox{\boldmath $\nu$}\times(\mbox{\boldmath $W$}_m\times\mbox{\boldmath $\nu$})
-\mbox{\boldmath $\nu$}\times\mbox{\boldmath $V$}_e\\
\\
\displaystyle
=-\frac{1}{\lambda}\mbox{\boldmath $\nu$}\times\{(\mbox{\boldmath $W$}_m-\mbox{\boldmath $V$}_m)\times\mbox{\boldmath $\nu$}\}
-\mbox{\boldmath $\nu$}\times\mbox{\boldmath $V$}_e
-\frac{1}{\lambda}\mbox{\boldmath $\nu$}\times(\mbox{\boldmath $V$}_m\times\mbox{\boldmath $\nu$})\\
\\
\displaystyle
=-\frac{1}{\lambda}\mbox{\boldmath $\nu$}\times(\mbox{\boldmath $R$}_m\times\mbox{\boldmath $\nu$})
-\mbox{\boldmath $\nu$}\times\mbox{\boldmath $V$}_e
-\frac{1}{\lambda}\mbox{\boldmath $\nu$}\times(\mbox{\boldmath $V$}_m\times\mbox{\boldmath $\nu$}).
\end{array}
$$
Then, from the second equation on (3.9) together with this yields
$$\begin{array}{l}
\displaystyle
\,\,\,\,\,\,
\mbox{\boldmath $\nu$}\times\mbox{\boldmath $R$}_e\cdot\nabla\times\mbox{\boldmath $R$}_e\\
\\
\displaystyle
=\left\{\frac{1}{\lambda}\mbox{\boldmath $\nu$}\times(\mbox{\boldmath $R$}_m\times\mbox{\boldmath $\nu$})
+\mbox{\boldmath $\nu$}\times\mbox{\boldmath $V$}_e
+\frac{1}{\lambda}\mbox{\boldmath $\nu$}\times(\mbox{\boldmath $V$}_m\times\mbox{\boldmath $\nu$})\right\}
\cdot
(\tau\mu\mbox{\boldmath $R$}_m+e^{-\tau T}\mu\mbox{\boldmath $G$})\\
\\
\displaystyle
=\frac{\tau\mu}{\lambda}
\mbox{\boldmath $\nu$}\times(\mbox{\boldmath $R$}_m\times\mbox{\boldmath $\nu$})\cdot\mbox{\boldmath $R$}_m
+\tau\mu\left\{\mbox{\boldmath $\nu$}\times\mbox{\boldmath $V$}_e
+\frac{1}{\lambda}\mbox{\boldmath $\nu$}\times(\mbox{\boldmath $V$}_m\times\mbox{\boldmath $\nu$})\right\}
\cdot\mbox{\boldmath $R$}_m\\
\\
\displaystyle
\,\,\,
+e^{-\tau T}\mu\left\{\frac{1}{\lambda}\mbox{\boldmath $\nu$}\times(\mbox{\boldmath $R$}_m\times\mbox{\boldmath $\nu$})
+\mbox{\boldmath $\nu$}\times\mbox{\boldmath $V$}_e
+\frac{1}{\lambda}\mbox{\boldmath $\nu$}\times(\mbox{\boldmath $V$}_m\times\mbox{\boldmath $\nu$})\right\}
\cdot\mbox{\boldmath $G$}\\
\\
\displaystyle
=\frac{\tau\mu}{\lambda}\vert\mbox{\boldmath $R$}_m\times\mbox{\boldmath $\nu$}\vert^2
+\tau\mu
\left(
\mbox{\boldmath $\nu$}\times(\mbox{\boldmath $V$}_e\times\mbox{\boldmath $\nu$})
+\frac{1}{\lambda}\mbox{\boldmath $V$}_m\times\mbox{\boldmath $\nu$}\right)
\cdot\mbox{\boldmath $R$}_m\times\mbox{\boldmath $\nu$}\\
\\
\displaystyle
\,\,\,
+e^{-\tau T}\mu
\left(\frac{1}{\lambda}\mbox{\boldmath $R$}_m\times\mbox{\boldmath $\nu$}
+
\mbox{\boldmath $\nu$}\times(\mbox{\boldmath $V$}_e\times\mbox{\boldmath $\nu$})
+\frac{1}{\lambda}
\mbox{\boldmath $V$}_m\times\mbox{\boldmath $\nu$}\right)\cdot\mbox{\boldmath $G$}\times\mbox{\boldmath $\nu$}\\
\\
\displaystyle
=\frac{\tau\mu}{\lambda}\vert\mbox{\boldmath $R$}_m\times\mbox{\boldmath $\nu$}\vert^2
+\mu
\left(\tau\mbox{\boldmath $V$}_{em}
+\frac{e^{-\tau T}}{\lambda}\mbox{\boldmath $G$}\times\mbox{\boldmath $\nu$}\right)
\cdot\mbox{\boldmath $R$}_m\times\mbox{\boldmath $\nu$}\\
\\
\displaystyle
\,\,\,
+e^{-\tau T}\mu\mbox{\boldmath $V$}_{em}\cdot\mbox{\boldmath $G$}\times\mbox{\boldmath $\nu$}.
\end{array}
$$
Thus we obtain
$$
\begin{array}{l}
\displaystyle
\,\,\,\,\,\,
\frac{1}{\mu\epsilon}
\int_{\partial D}\mbox{\boldmath $\nu$}\times\mbox{\boldmath $R$}_e\cdot\nabla\times\mbox{\boldmath $R$}_edS
\\
\\
\displaystyle
=\frac{\tau}{\epsilon}
\int_{\partial D}\frac{1}{\lambda}\vert\mbox{\boldmath $R$}_m\times\mbox{\boldmath $\nu$}\vert^2dS
+\frac{\tau}{\epsilon}\int_{\partial D}
\left(
\mbox{\boldmath $V$}_{em}
+\frac{e^{-\tau T}}{\tau\lambda}\mbox{\boldmath $G$}\times\mbox{\boldmath $\nu$}\right)
\cdot\mbox{\boldmath $R$}_m\times\mbox{\boldmath $\nu$}dS\\
\\
\displaystyle
\,\,\,
+\frac{e^{-\tau T}}{\epsilon}
\int_{\partial D}
\mbox{\boldmath $V$}_{em}
\cdot\mbox{\boldmath $G$}\times\mbox{\boldmath $\nu$}dS
\\
\\
\displaystyle
=\frac{\tau}{\epsilon}\int_{\partial D}\frac{1}{\lambda}
\left\vert
\mbox{\boldmath $R$}_m\times\mbox{\boldmath $\nu$}
+\frac{\lambda}{2}
\left(\mbox{\boldmath $V$}_{em}+\frac{e^{-\tau T}}{\tau\lambda}\mbox{\boldmath $G$}\times\mbox{\boldmath $\nu$}\right)\right\vert^2
dS\\
\\
\displaystyle
\,\,\,
-\frac{\tau}{4\epsilon}\int_{\partial D}\lambda
\left\vert\mbox{\boldmath $V$}_{em}
+\frac{e^{-\tau T}}{2\lambda}\mbox{\boldmath $G$}\times\mbox{\boldmath $\nu$}\right\vert^2dS
+\frac{e^{-\tau T}}{\epsilon}
\int_{\partial D}\mbox{\boldmath $V$}_{em}\cdot\mbox{\boldmath $G$}\times\mbox{\boldmath $\nu$}dS\\
\\
\displaystyle
\ge
\frac{\tau}{2\epsilon}\int_{\partial D}
\frac{1}{\lambda}
\vert\mbox{\boldmath $R$}_m\times\mbox{\boldmath $\nu$}\vert^2 dS
-\frac{\tau}{2\epsilon}
\int_{\partial D}\lambda
\left\vert
\mbox{\boldmath $V$}_{em}
+\frac{e^{-\tau T}}{\tau\lambda}
\mbox{\boldmath $G$}\times\mbox{\boldmath $\nu$}\right\vert^2dS\\
\\
\displaystyle
\,\,\,
+\frac{e^{-\tau T}}{\epsilon}
\int_{\partial D}\mbox{\boldmath $V$}_{em}\cdot\mbox{\boldmath $G$}\times\mbox{\boldmath $\nu$}dS.
\end{array}
\tag {3.52}
$$
Now a combination of (3.51) and (3.52) gives
$$\begin{array}{ll}
\displaystyle
\frac{1}{2}\tilde{E}_e(\tau) & \displaystyle
\le\frac{\tau}{2\epsilon}
\int_{\partial D}\lambda
\left\vert
\mbox{\boldmath $V$}_{em}
+\frac{e^{-\tau T}}{\tau\lambda}
\mbox{\boldmath $G$}\times\mbox{\boldmath $\nu$}\right\vert^2dS\\
\\
\displaystyle
 & \,\,\,\displaystyle
 +
\frac{e^{-2\tau T}}{2\tau^2}\int_{\Bbb R^3\setminus\overline D}\vert\mbox{\boldmath $F$}_e-\mbox{\boldmath $F$}_e^0\vert^2dx
-\frac{e^{-\tau T}}{\epsilon}
\int_{\partial D}\mbox{\boldmath $V$}_{em}\cdot\mbox{\boldmath $G$}\times\mbox{\boldmath $\nu$}dS
\end{array}
$$
and hence
$$
\begin{array}{ll}
\displaystyle
\tilde{E}_e(\tau)
&
\displaystyle
\le \frac{\tau}{\epsilon}
\int_{\partial D}\lambda
\left\vert
\mbox{\boldmath $V$}_{em}
+\frac{e^{-\tau T}}{\tau\lambda}
\mbox{\boldmath $G$}\times\mbox{\boldmath $\nu$}\right\vert^2dS\\
\\
\displaystyle
& \,\,\,\displaystyle
+
\frac{e^{-2\tau T}}{\tau^2}\int_{\Bbb R^3\setminus\overline D}\vert\mbox{\boldmath $F$}_e-\mbox{\boldmath $F$}_e^0\vert^2dx
-\frac{2e^{-\tau T}}{\epsilon}
\int_{\partial D}\mbox{\boldmath $V$}_{em}\cdot\mbox{\boldmath $G$}\times\mbox{\boldmath $\nu$}dS\\
\\
\displaystyle
&
= \displaystyle
\frac{\tau}{\epsilon}
\int_{\partial D}\lambda
\left(\vert
\mbox{\boldmath $V$}_{em}\vert^2
+\left\vert\frac{e^{-\tau T}}{\tau\lambda}
\mbox{\boldmath $G$}\times\mbox{\boldmath $\nu$}\right\vert^2\right)dS
+\frac{e^{-2\tau T}}{\tau^2}\int_{\Bbb R^3\setminus\overline D}\vert\mbox{\boldmath $F$}_e-\mbox{\boldmath $F$}_e^0\vert^2dx.
\end{array}
$$
Then this together with (3.8) and (3.10) yields (3.49).

\section{Proof of Lemmas 2.2 and 2.3}

First we describe the proof of Lemma 2.2.

\subsection{A reduction}

\proclaim{\noindent Lemma 4.1}
We have
$$
\begin{array}{ll}
\displaystyle
\tilde{J}_e(\tau)
&
\displaystyle
=\frac{1}{\mu\epsilon}
\int_{\partial D}\frac{1}{c}
\left\{c\,\mbox{\boldmath $\nu$}\times(\mbox{\boldmath $V$}_e^0\times\mbox{\boldmath $\nu$})
-(\nabla\times\mbox{\boldmath $V$}_e^0)\times\mbox{\boldmath $\nu$}\right\}
\cdot(\nabla\times\mbox{\boldmath $V$}_e^0)\times\mbox{\boldmath $\nu$}\,dS
\\
\\
\displaystyle
&
\,\,\,\displaystyle
+O(\tau^{-1/2}e^{-\tau T})
\end{array}
\tag {4.1}
$$
and
$$
\begin{array}{l}
\displaystyle
\,\,\,\,\,\,
\tilde{J}_{e}(\tau)
+\frac{1}{\mu\epsilon}
\int_{\partial D}c\vert\mbox{\boldmath $V$}_{em}\vert^2\,dS\\
\\
\displaystyle
=\frac{1}{\mu\epsilon}
\int_{\partial D}
\left\{
c\mbox{\boldmath $\nu$}\times(\mbox{\boldmath $V$}_e^0\times\mbox{\boldmath $\nu$})
-(\nabla\times\mbox{\boldmath $V$}_e^0)\times\mbox{\boldmath $\nu$}\right\}
\cdot
\mbox{\boldmath $\nu$}\times(\mbox{\boldmath $V$}_e^0\times\mbox{\boldmath $\nu$})\,dS
+O(\tau^{-1/2}e^{-\tau T}),
\end{array}
\tag {4.2}
$$
where $c=c(x,\tau)\equiv \tau\mu\lambda(x)$ and
$\mbox{\boldmath $V$}_{em}$ is given by (3.47).

\endproclaim

{\it\noindent Proof.}
We prepare two elementary estimates.
Since $\mbox{\boldmath $V$}_e^0$ satisfies (1.14), we have
$$\displaystyle
\frac{1}{\mu\epsilon}
\nabla\times\nabla\times\mbox{\boldmath $V$}_e^0
+\tau^2\mbox{\boldmath $V$}_e^0=\mbox{\boldmath $0$}\,\,\text{in}\,D.
$$
Thus (3.36) gives
$$\displaystyle
\Vert\nabla\times\nabla\times\mbox{\boldmath $V$}_e^0\Vert_{L^2(D)}
=O(\tau^{-1/2}).
$$
By the trace theorem \cite{DL3} this together with (3.37) yields
$$\displaystyle
\Vert(\nabla\times\mbox{\boldmath $V$}_e^0)\times\mbox{\boldmath $\nu$}\Vert_{L^2(\partial D)}
=O(\tau^{-1/2}).
\tag {4.3}
$$
Next from  equation (3.38) and (3.39) we have
$$\displaystyle
\Vert\nabla\times\nabla\times(\mbox{\boldmath $V$}_e-\mbox{\boldmath $V$}_e^0)\Vert_{L^2(\Bbb R^3)}
=O(\tau e^{-\tau T}).
$$
Again by the trace theorem \cite{DL3} this together with (3.40) gives
$$\displaystyle
\Vert(\nabla\times(\mbox{\boldmath $V$}_e-\mbox{\boldmath $V$}_e^0))\times\mbox{\boldmath $\nu$}\Vert_{L^2(\partial D)}
=O(\tau e^{-\tau T}).
\tag {4.4}
$$

Moreover, we add also
$$\displaystyle
\Vert\mbox{\boldmath $V$}_e^0\times\mbox{\boldmath $\nu$}\Vert_{L^2(\partial D)}
=O(\tau^{-3/2})
\tag {4.5}
$$
and
$$\displaystyle
\Vert
\mbox{\boldmath $\nu$}\times(\mbox{\boldmath $V$}_e-\mbox{\boldmath $V$}_e^0)
\Vert_{L^2(\partial D)}
=O(e^{-\tau T}),
\tag {4.6}
$$
which are derived from the trace theorem \cite{DL3} together
with  (3.36) and (3.37); (3.39) and (3.40), respectively.

Then, from (4.3) and (4.4) we have, in $L^2(\partial D)$
$$\begin{array}{l}
\displaystyle
\,\,\,\,\,\,
\mbox{\boldmath $V$}_m\times\mbox{\boldmath $\nu$}
\\
\\
\displaystyle
=-\frac{1}{\tau\mu}(\nabla\times\mbox{\boldmath $V$}_e)\times\mbox{\boldmath $\nu$}
+\frac{e^{-\tau T}}{\tau}\mbox{\boldmath $H$}_0\times\mbox{\boldmath $\nu$}\\
\\
\displaystyle
=-\frac{1}{\tau\mu}(\nabla\times\mbox{\boldmath $V$}_e^0)\times\mbox{\boldmath $\nu$}
-\frac{1}{\tau\mu}\{\nabla\times(\mbox{\boldmath $V$}_e-\mbox{\boldmath $V$}_e^0)\}\times\mbox{\boldmath $\nu$}
+O(\tau^{-1}e^{-\tau T})\\
\\
\displaystyle
=-\frac{1}{\tau\mu}(\nabla\times\mbox{\boldmath $V$}_e^0)\times\mbox{\boldmath $\nu$}
+O(e^{-\tau T}).
\end{array}
$$
Thus we have
$$\begin{array}{ll}
\displaystyle
\int_{\partial D}\frac{\tau\mu}{\lambda}\vert\mbox{\boldmath $V$}_m\times\mbox{\boldmath $\nu$}\vert^2dS
& \displaystyle
=\int_{\partial D}\frac{\tau\mu}{\lambda}\frac{1}{\tau^2\mu^2}
\vert(\nabla\times\mbox{\boldmath $V$}_e^0)\times\mbox{\boldmath $\nu$}\vert^2 dS
\\
\\
\displaystyle
& \,\,\,\displaystyle
+O(\tau^{-1/2}e^{-\tau T})+O(e^{-2\tau T})\\
\\
\displaystyle
&
\displaystyle
=\int_{\partial D}\frac{1}{c}\vert(\nabla\times\mbox{\boldmath $V$}_e^0)\times\mbox{\boldmath $\nu$}\vert^2 dS
+O(\tau^{-1/2}e^{-\tau T}).
\end{array}
\tag {4.7}
$$

Second write
$$\begin{array}{l}
\displaystyle
\,\,\,\,\,\,
(\mbox{\boldmath $\nu$}\times\mbox{\boldmath $V$}_e)\cdot\nabla\times\mbox{\boldmath $V$}_e
\\
\\
\displaystyle
=
(\mbox{\boldmath $\nu$}\times\mbox{\boldmath $V$}_e)\cdot\mbox{\boldmath $\nu$}\times\{(\nabla\times\mbox{\boldmath $V$}_e)\times\mbox{\boldmath $\nu$}\}
\\
\\
\displaystyle
=(\mbox{\boldmath $\nu$}\times\mbox{\boldmath $V$}_e^0)
\cdot\mbox{\boldmath $\nu$}\times\{(\nabla\times\mbox{\boldmath $V$}_e)\times\mbox{\boldmath $\nu$}\}
+
\{\mbox{\boldmath $\nu$}\times(\mbox{\boldmath $V$}_e-\mbox{\boldmath $V$}_e^0)\}\cdot\mbox{\boldmath $\nu$}\times\{(\nabla\times\mbox{\boldmath $V$}_e)\times\mbox{\boldmath $\nu$}\}
\\
\\
\displaystyle
=(\mbox{\boldmath $\nu$}\times\mbox{\boldmath $V$}_e^0)
\cdot\mbox{\boldmath $\nu$}\times\{(\nabla\times\mbox{\boldmath $V$}_e^0)\times\mbox{\boldmath $\nu$}\}
+
(\mbox{\boldmath $\nu$}\times\mbox{\boldmath $V$}_e^0)
\cdot\mbox{\boldmath $\nu$}\times\{(\nabla\times
(\mbox{\boldmath $V$}_e-\mbox{\boldmath $V$}_e^0))\times\mbox{\boldmath $\nu$}\}\\
\\
\displaystyle
\,\,\,
+\{\mbox{\boldmath $\nu$}\times(\mbox{\boldmath $V$}_e-\mbox{\boldmath $V$}_e^0)\}\cdot\mbox{\boldmath $\nu$}\times\{(\nabla\times\mbox{\boldmath $V$}_e^0)\times\mbox{\boldmath $\nu$}\}
\\
\\
\displaystyle
\,\,\,
+\{\mbox{\boldmath $\nu$}\times(\mbox{\boldmath $V$}_e-\mbox{\boldmath $V$}_e^0)\}\cdot\mbox{\boldmath $\nu$}\times\{(\nabla\times(\mbox{\boldmath $V$}_e-\mbox{\boldmath $V$}_e^0))\times\mbox{\boldmath $\nu$}\}.
\end{array}
$$
Then using (4.3), (4.4), (4.5) and (4.6), we have
$$\begin{array}{l}
\displaystyle
\,\,\,\,\,\,
\int_{\partial D}(\mbox{\boldmath $\nu$}\times\mbox{\boldmath $V$}_e)\cdot\nabla\times\mbox{\boldmath $V$}_e
dS
\\
\\
\displaystyle
=\int_{\partial D}(\mbox{\boldmath $\nu$}\times\mbox{\boldmath $V$}_e^0)
\cdot\mbox{\boldmath $\nu$}\times\{(\nabla\times\mbox{\boldmath $V$}_e^0)\times\mbox{\boldmath $\nu$}\}
dS\\
\\
\displaystyle
\,\,\,
+O(\tau^{-3/2}\tau e^{-\tau T})
+O(e^{-\tau T}\tau^{-1/2})
+O(e^{-\tau T}\tau e^{-\tau T})
=O(\tau^{-1/2}e^{-\tau T}).
\end{array}
\tag {4.8}
$$
Since
$$\displaystyle
(\mbox{\boldmath $\nu$}\times\mbox{\boldmath $V$}_e^0)
\cdot\mbox{\boldmath $\nu$}\times\{(\nabla\times\mbox{\boldmath $V$}_e^0)\times\mbox{\boldmath $\nu$}\}
=\mbox{\boldmath $\nu$}\times(\mbox{\boldmath $V$}_e^0\times\mbox{\boldmath $\nu$})
\cdot
\{(\nabla\times\mbox{\boldmath $V$}_e^0)\times\mbox{\boldmath $\nu$}\}
$$
and
$$\displaystyle
\tilde{J}_e(\tau)
=\frac{1}{\mu\epsilon}
\left
(\int_{\partial D}(\mbox{\boldmath $\nu$}\times\mbox{\boldmath $V$}_e)\cdot\nabla\times\mbox{\boldmath $V$}_e
-\frac{\tau\mu}{\lambda}
\vert\mbox{\boldmath $V$}_m\times\mbox{\boldmath $\nu$}\vert^2\right)dS,
$$
from (4.7) and (4.8) we obtain
$$\begin{array}{ll}
\displaystyle
\tilde{J}_e(\tau)
&
\displaystyle
=\frac{1}{\mu\epsilon}
\int_{\partial D}
\left\{\mbox{\boldmath $\nu$}\times
(\mbox{\boldmath $V$}_e^0\times\mbox{\boldmath $\nu$})
\cdot\{(\nabla\times\mbox{\boldmath $V$}_e^0)\times\mbox{\boldmath $\nu$}\}
-\frac{1}{c}\vert(\nabla\times\mbox{\boldmath $V$}_e^0)\times\mbox{\boldmath $\nu$}\vert^2 \right\}dS\\
\\
\displaystyle
&
\displaystyle
\,\,\,
+O(\tau^{-1/2}e^{-\tau T}).
\end{array}
$$
This is nothing but (4.1).

Next, in $L^2(\partial D)$ from (4.3), (4.4) and (4.6)  we have
$$\begin{array}{l}
\displaystyle
\,\,\,\,\,\,
\mbox{\boldmath $V$}_{em}
\\
\\
\displaystyle
=\mbox{\boldmath $\nu$}\times(\mbox{\boldmath $V$}_e\times\mbox{\boldmath $\nu$})
+\frac{1}{\lambda}\mbox{\boldmath $V$}_m\times\mbox{\boldmath $\nu$}\\
\\
\displaystyle
=\mbox{\boldmath $\nu$}\times(\mbox{\boldmath $V$}_e\times\mbox{\boldmath $\nu$})
+\frac{1}{\lambda}
\left\{
-\frac{1}{\tau\mu}
(\nabla\times\mbox{\boldmath $V$}_e)\times\mbox{\boldmath $\nu$}
-\frac{e^{-\tau T}}{\tau}\mbox{\boldmath $H$}_0\right\}
\times\mbox{\boldmath $\nu$}\\
\\
\displaystyle
=\mbox{\boldmath $\nu$}\times(\mbox{\boldmath $V$}_e\times\mbox{\boldmath $\nu$})
-\frac{1}{c}(\nabla\times\mbox{\boldmath $V$}_e)\times\mbox{\boldmath $\nu$}
-\frac{e^{-\tau T}}{\lambda\tau}\mbox{\boldmath $H$}_0\times\mbox{\boldmath $\nu$}\\
\\
\displaystyle
=\mbox{\boldmath $\nu$}\times(\mbox{\boldmath $V$}_e^0\times\mbox{\boldmath $\nu$})
+\mbox{\boldmath $\nu$}\times\{(\mbox{\boldmath $V$}_e-\mbox{\boldmath $V$}_e^0)\times\mbox{\boldmath $\nu$}\}
\\
\\
\displaystyle
\,\,\,
-\frac{1}{c}(\nabla\times\mbox{\boldmath $V$}_e^0)\times\mbox{\boldmath $\nu$}
-\frac{1}{c}\{\nabla\times(\mbox{\boldmath $V$}_e-\mbox{\boldmath $V$}_e^0)\}\times\mbox{\boldmath $\nu$}
+O(\tau^{-1}e^{-\tau T})
\\
\\
\displaystyle
=\mbox{\boldmath $\nu$}\times(\mbox{\boldmath $V$}_e^0\times\mbox{\boldmath $\nu$})
+O(e^{-\tau T})
-\frac{1}{c}(\nabla\times\mbox{\boldmath $V$}_e^0)\times\mbox{\boldmath $\nu$}+O(e^{-\tau T})
+O(\tau^{-1}e^{-\tau T})\\
\\
\displaystyle
=\mbox{\boldmath $\nu$}\times(\mbox{\boldmath $V$}_e^0\times\mbox{\boldmath $\nu$})
-\frac{1}{c}(\nabla\times\mbox{\boldmath $V$}_e^0)\times\mbox{\boldmath $\nu$}
+O(e^{-\tau T}).
\end{array}
$$
Then from (4.3) and (4.5)  we obtain
$$\begin{array}{l}
\,\,\,\,\,\,
\displaystyle
\frac{1}{\mu\epsilon}
\int_{\partial D}c\vert\mbox{\boldmath $V$}_{em}\vert^2 dS
\\
\\
\displaystyle
=\frac{1}{\mu\epsilon}\int_{\partial D}
c\left\vert
\mbox{\boldmath $\nu$}\times(\mbox{\boldmath $V$}_e^0\times\mbox{\boldmath $\nu$})
-\frac{1}{c}(\nabla\times\mbox{\boldmath $V$}_e^0)\times\mbox{\boldmath $\nu$}
\right\vert^2 dS\\
\\
\displaystyle
\,\,\,
+O(\tau\tau^{-3/2}e^{-\tau T})
+O(\tau e^{-2\tau T})\\
\\
\displaystyle
=\frac{1}{\mu\epsilon}\int_{\partial D}
\frac{1}{c}\left\vert
c\mbox{\boldmath $\nu$}\times(\mbox{\boldmath $V$}_e^0\times\mbox{\boldmath $\nu$})
-(\nabla\times\mbox{\boldmath $V$}_e^0)\times\mbox{\boldmath $\nu$}
\right\vert^2 dS
+O(\tau^{-1/2}e^{-\tau T}).
\end{array}
\tag {4.9}
$$
Now a combination of (4.1) and (4.9) yields (4.2).

\subsection{Preliminary computation for the proof of (i) and (ii)}

\subsubsection{Computation of $\mbox{\boldmath $\nu$}\times(\mbox{\boldmath$V$}_e^0\times\mbox{\boldmath $\nu$})$}

By (18) in \cite{IEM} we know that $\mbox{\boldmath $V$}_e^0$ is smooth outside $B$ and has the form
$$\begin{array}{ll}
\displaystyle
\mbox{\boldmath $V$}_e^0
=K(\tau)\tilde{f}(\tau)v\mbox{\boldmath $M$}\mbox{\boldmath $a$}
& \text{in}\,\Bbb R^3\setminus\overline B,
\end{array}
\tag {4.10}
$$
where
$$
\left\{
\begin{array}{l}
\displaystyle
v=v(x)=\frac{e^{-\tau\sqrt{\mu\epsilon}\vert x-p\vert}}{\vert x-p\vert},
\\
\\
\displaystyle
K(\tau)=\frac{\mu\tau\varphi(\tau\sqrt{\mu\epsilon})}{(\tau\sqrt{\mu\epsilon})^3},
\\
\\
\displaystyle
\varphi(\xi)=\xi\cosh\xi-\sinh\xi
\end{array}
\right.
$$
and
$$
\left\{
\begin{array}{l}
\displaystyle
\mbox{\boldmath $M$}
=\mbox{\boldmath $M$}(x;\tau)
=AI_3-B\,\frac{x-p}{\vert x-p\vert}\otimes
\frac{x-p}{\vert x-p\vert},
\\
\\
\displaystyle
A=A(x,\tau)=1+\frac{1}{\tau\sqrt{\mu\epsilon}}
\left(\frac{1}{\vert x-p\vert}
+\frac{1}{\tau\sqrt{\mu\epsilon}\vert x-p\vert^2}\right),\\
\\
\displaystyle
B=B(x,\tau)=1+\frac{3}{\tau\sqrt{\mu\epsilon}}
\left(\frac{1}{\vert x-p\vert}
+\frac{1}{\tau\sqrt{\mu\epsilon}\vert x-p\vert^2}\right).
\end{array}
\right.
$$
The expression (4.10) is a simple application of the mean value theorem \cite{CH} for the modified
Helmholtz equation and a special form of the fundamental solution
of equation (1.14).

From (4.10) we have
$$\displaystyle
\mbox{\boldmath $\nu$}\times(\mbox{\boldmath$V$}_e^0\times\mbox{\boldmath $\nu$})
=K(\tau)\tilde{f}(\tau)v\,\mbox{\boldmath $\nu$}\times\{(\mbox{\boldmath $M$}\mbox{\boldmath $a$})\times\mbox{\boldmath $\nu$}\}.
$$

\subsubsection{Computation of $(\nabla\times\mbox{\boldmath $V$}_e^0)\times\mbox{\boldmath $\nu$}$}

By (23) in \cite{IEM}, from (4.10) we have already derived
the expression of $\nabla\times\mbox{\boldmath $V$}_e^0$ outside $B$:
$$\displaystyle
\nabla\times\mbox{\boldmath $V$}_e^0
=-\tau\sqrt{\mu\epsilon}\,K(\tau)\tilde{f}(\tau)v
\,\left(1+\frac{1}{\tau\sqrt{\mu\epsilon}\vert x-p\vert}\right)
\frac{x-p}{\vert x-p\vert}\times\mbox{\boldmath $a$}.
$$
Thus we obtain
$$\displaystyle
(\nabla\times\mbox{\boldmath $V$}_e^0)\times\mbox{\boldmath $\nu$}
=-\tau\sqrt{\mu\epsilon}\,K(\tau)\tilde{f}(\tau)v
\,\left(1+\frac{1}{\tau\sqrt{\mu\epsilon}\vert x-p\vert}\right)
\left(\frac{x-p}{\vert x-p\vert}\times\mbox{\boldmath $a$}\right)\times\mbox{\boldmath $\nu$}.
$$

\subsubsection{Computation of $\mbox{\boldmath $\nu$}\times\{(\mbox{\boldmath $M$}\mbox{\boldmath $a$})\times\mbox{\boldmath $\nu$}\}$}

We have
$$\displaystyle
\mbox{\boldmath $M$}\mbox{\boldmath $a$}
=A\,\mbox{\boldmath $a$}-B\,\left(\mbox{\boldmath $a$}\cdot\frac{x-p}{\vert x-p\vert}\right)\,\frac{x-p}{\vert x-p\vert}
$$
and
$$\displaystyle
(\mbox{\boldmath $M$}\mbox{\boldmath $a$})\cdot\mbox{\boldmath $\nu$}
=A\,\mbox{\boldmath $a$}\cdot\mbox{\boldmath $\nu$}
-B\left(\mbox{\boldmath $a$}\cdot\frac{x-p}{\vert x-p\vert}\right)
\left(\mbox{\boldmath $\nu$}\cdot\frac{x-p}{\vert x-p\vert}\right).
$$
Combining these with expression
$$\begin{array}{c}
\displaystyle
\mbox{\boldmath $\nu$}\times\{(\mbox{\boldmath $M$}\mbox{\boldmath $a$})\times\mbox{\boldmath $\nu$}\}
=\mbox{\boldmath $M$}\mbox{\boldmath $a$}-\{(\mbox{\boldmath $M$}\mbox{\boldmath $a$})\cdot\mbox{\boldmath $\nu$}\}\mbox{\boldmath $\nu$},
\end{array}
$$
we obtain
$$\displaystyle
\mbox{\boldmath $\nu$}\times\{(\mbox{\boldmath $M$}\mbox{\boldmath $a$})\times\mbox{\boldmath $\nu$}\}
=A\,\mbox{\boldmath $\nu$}\times(\mbox{\boldmath $a$}\times\mbox{\boldmath $\nu$})
-B\left(\mbox{\boldmath $a$}\cdot\frac{x-p}{\vert x-p\vert}\right)
\mbox{\boldmath $\nu$}\times
\left(\frac{x-p}{\vert x-p\vert}\times\mbox{\boldmath $\nu$}\right).
$$

\subsubsection{Computation of $\left(\frac{x-p}{\vert x-p\vert}\times\mbox{\boldmath $a$}\right)\times\mbox{\boldmath $\nu$}$}

From the identity $\mbox{\boldmath $A$}\times(\mbox{\boldmath $B$}\times\mbox{\boldmath $C$})
=(\mbox{\boldmath $A$}\cdot\mbox{\boldmath $C$})\mbox{\boldmath $B$}
-(\mbox{\boldmath $A$}\cdot\mbox{\boldmath $B$})\mbox{\boldmath $C$}$,
we obtain
$$\begin{array}{l}
\displaystyle
\,\,\,\,\,\,
\left(\frac{x-p}{\vert x-p\vert}\times\mbox{\boldmath $a$}\right)\times\mbox{\boldmath $\nu$}
\\
\\
\displaystyle
=-\mbox{\boldmath $\nu$}\times\left(\frac{x-p}{\vert x-p\vert}
\times\mbox{\boldmath $a$}\right)\\
\\
\displaystyle
=\mbox{\boldmath $\nu$}\times
\left(\mbox{\boldmath $a$}\times\frac{x-p}{\vert x-p\vert}\right)
\\
\\
\displaystyle
=\left(\mbox{\boldmath $\nu$}\cdot\frac{x-p}{\vert x-p\vert}\right)\mbox{\boldmath $a$}
-\left(\mbox{\boldmath $\nu$}\cdot\mbox{\boldmath $a$}\right)
\frac{x-p}{\vert x-p\vert}
\\
\\
\displaystyle
=\left(\mbox{\boldmath $\nu$}\cdot\frac{x-p}{\vert x-p\vert}\right)
\mbox{\boldmath $\nu$}\times(\mbox{\boldmath $a$}\times\mbox{\boldmath $\nu$})
-\left(\mbox{\boldmath $\nu$}\cdot\mbox{\boldmath $a$}\right)
\mbox{\boldmath $\nu$}\times\left(\frac{x-p}{\vert x-p\vert}\times\mbox{\boldmath $\nu$}\right)
\end{array}
$$

\subsubsection{Computation of
$
c\,\mbox{\boldmath $\nu$}\times(\mbox{\boldmath$V$}_e^0\times\mbox{\boldmath $\nu$})
-(\nabla\times\mbox{\boldmath $V$}_e^0)\times\mbox{\boldmath $\nu$}
$}

We have
$$\begin{array}{l}
\displaystyle
\,\,\,\,\,\,
c\,\mbox{\boldmath $\nu$}\times(\mbox{\boldmath$V$}_e^0\times\mbox{\boldmath $\nu$})
-(\nabla\times\mbox{\boldmath $V$}_e^0)\times\mbox{\boldmath $\nu$}\\
\\
\displaystyle
=K(\tau)\tilde{f}(\tau)v\\
\\
\displaystyle
\,\,\,
\times
\left\{
c\,\mbox{\boldmath $\nu$}\times\{(\mbox{\boldmath $M$}\mbox{\boldmath $a$})\times\mbox{\boldmath $\nu$}\}
+\tau\sqrt{\mu\epsilon}
\left(1+\frac{1}{\tau\sqrt{\mu\epsilon}\vert x-p\vert}\right)
\left(\frac{x-p}{\vert x-p\vert}\times\mbox{\boldmath $a$}\right)\times\mbox{\boldmath $\nu$}\right\}\\
\\
\displaystyle
=K(\tau)\tilde{f}(\tau)\tau\mu\,\\
\\
\displaystyle
\,\,\,
\times
\left\{
\lambda(x)\,\mbox{\boldmath $\nu$}\times\{(\mbox{\boldmath $M$}\mbox{\boldmath $a$})\times\mbox{\boldmath $\nu$}\}
+\sqrt{\frac{\epsilon}{\mu}}\,
\left(1+\frac{1}{\tau\sqrt{\mu\epsilon}\vert x-p\vert}\right)
\left(\frac{x-p}{\vert x-p\vert}\times\mbox{\boldmath $a$}\right)\times\mbox{\boldmath $\nu$}\right\}\\
\\
\displaystyle
=K(\tau)\tilde{f}(\tau)\tau\mu\,v\,
\left(
P\,\mbox{\boldmath $X$}
-Q\,\mbox{\boldmath $Y$}
\right),
\end{array}
$$
where
$$
\left\{
\begin{array}{l}
\displaystyle
P=\lambda\,A+\sqrt{\frac{\epsilon}{\mu}}\left(1+\frac{1}{\tau\sqrt{\mu\epsilon}\vert x-p\vert}\right)
\left(\mbox{\boldmath $\nu$}\cdot\frac{x-p}{\vert x-p\vert}\right),\\
\\
\displaystyle
Q=\lambda\,B\left(\mbox{\boldmath $a$}\cdot\frac{x-p}{\vert x-p\vert}\right)
+\sqrt{\frac{\epsilon}{\mu}}
\left(1+\frac{1}{\tau\sqrt{\mu\epsilon}\vert x-p\vert}\right)
(\mbox{\boldmath $\nu$}\cdot\mbox{\boldmath $a$})
\end{array}
\right.
$$
and
$$\begin{array}{c}
\displaystyle
\mbox{\boldmath $X$}=\mbox{\boldmath $\nu$}\times(\mbox{\boldmath $a$}\times\mbox{\boldmath $\nu$}),\,\,
\mbox{\boldmath $Y$}=\mbox{\boldmath $\nu$}\times\left(\frac{x-p}{\vert x-p\vert}\times\mbox{\boldmath $\nu$}\right).
\end{array}
$$

\subsubsection{Computation of
$
\left\{c\,\mbox{\boldmath $\nu$}\times(\mbox{\boldmath$V$}_e^0\times\mbox{\boldmath $\nu$})
-(\nabla\times\mbox{\boldmath $V$}_e^0)\times\mbox{\boldmath $\nu$}\right\}
\cdot(\nabla\times\mbox{\boldmath $V$}_e^0)\times\mbox{\boldmath $\nu$}
$
}

We have
$$\displaystyle
(\nabla\times\mbox{\boldmath $V$}_e^0)\times\mbox{\boldmath $\nu$}
=-\tau\sqrt{\mu\epsilon}
K(\tau)\tilde{f}(\tau)v(x)\,
\left(1+\frac{1}{\tau\sqrt{\mu\epsilon}\vert x-p\vert}\right)
\,
\left\{
\left(\mbox{\boldmath $\nu$}\cdot\frac{x-p}{\vert x-p\vert}\right)\mbox{\boldmath $X$}
-(\mbox{\boldmath $\nu$}\cdot\mbox{\boldmath $a$})\mbox{\boldmath $Y$}\right\}.
$$
Thus we obtain
$$
\begin{array}{l}
\displaystyle
\,\,\,\,\,\,
\left\{c\,\mbox{\boldmath $\nu$}\times(\mbox{\boldmath$V$}_e^0\times\mbox{\boldmath $\nu$})
-(\nabla\times\mbox{\boldmath $V$}_e^0)\times\mbox{\boldmath $\nu$}\right\}
\cdot(\nabla\times\mbox{\boldmath $V$}_e^0)\times\mbox{\boldmath $\nu$}\\
\\
\displaystyle
=-\tau^2\mu\sqrt{\mu\epsilon}K(\tau)^2(\tilde{f}(\tau))^2v(x)^2\,\left(1+\frac{1}{\tau\sqrt{\mu\epsilon}\vert x-p\vert}\right)
\\
\\
\displaystyle
\,\,\,
\times
\left(
P\,\mbox{\boldmath $X$}
-Q\,\mbox{\boldmath $Y$}
\right)\cdot\left\{
\left(\mbox{\boldmath $\nu$}\cdot\frac{x-p}{\vert x-p\vert}\right)\mbox{\boldmath $X$}
-(\mbox{\boldmath $\nu$}\cdot\mbox{\boldmath $a$})\mbox{\boldmath $Y$}\right\}.
\end{array}
\tag {4.11}
$$

\subsubsection{Computation of
$
\left\{c\,\mbox{\boldmath $\nu$}\times(\mbox{\boldmath$V$}_e^0\times\mbox{\boldmath $\nu$})
-(\nabla\times\mbox{\boldmath $V$}_e^0)\times\mbox{\boldmath $\nu$}\right\}
\cdot
\mbox{\boldmath $\nu$}\times(\mbox{\boldmath $V$}_e^0\times\mbox{\boldmath $\nu$})
$
}

We have
$$\displaystyle
\mbox{\boldmath $\nu$}\times(\mbox{\boldmath $V$}_e^0\times\mbox{\boldmath $\nu$})
=K(\tau)\tilde{f}(\tau) v(x)\,
\left\{
A\,\mbox{\boldmath $X$}
-B\,
\left(\mbox{\boldmath $a$}\cdot\frac{x-p}{\vert x-p\vert}\right)
\mbox{\boldmath $Y$}
\right\}.
$$
Therefore we obtain
$$
\begin{array}{l}
\displaystyle
\,\,\,\,\,\,
\left\{c\,\mbox{\boldmath $\nu$}\times(\mbox{\boldmath$V$}_e^0\times\mbox{\boldmath $\nu$})
-(\nabla\times\mbox{\boldmath $V$}_e^0)\times\mbox{\boldmath $\nu$}\right\}
\cdot
\mbox{\boldmath $\nu$}\times(\mbox{\boldmath $V$}_e^0\times\mbox{\boldmath $\nu$})\\
\\
\displaystyle
=\tau\mu\,K(\tau)^2(\tilde{f}(\tau))^2v(x)^2\,
\left(
P\,\mbox{\boldmath $X$}
-Q\,\mbox{\boldmath $Y$}
\right)
\cdot
\left\{
A\,\mbox{\boldmath $X$}
-B\,
\left(\mbox{\boldmath $a$}\cdot\frac{x-p}{\vert x-p\vert}\right)
\mbox{\boldmath $Y$}
\right\}.
\end{array}
\tag {4.12}
$$

\subsection{Proof of (i)}

In this and next subsections,
to emphasize the dependence on space variable $x$
we write $P=P_x$, $Q=Q_x$, $\mbox{\boldmath $X$}=\mbox{\boldmath $X$}_x$,
$\mbox{\boldmath $Y$}=\mbox{\boldmath $Y$}_x$ and $\mbox{\boldmath $\nu$}=\mbox{\boldmath $\nu$}_x$.

First we prepare an elementary lemma which can be proved by using a contradiction argument.

\proclaim{\noindent Lemma 4.2}
Given $\eta>0$ there exists a positive number $\delta$ such that,
for all $x\in\partial D\cap B_{d_{\partial D}(p)+\delta}(p)$
$$\displaystyle
\mbox{\boldmath $\nu$}_x\cdot\frac{p-x}{\vert x-p\vert}>1-\eta.
$$
\endproclaim

Since
$$\displaystyle
\vert\mbox{\boldmath $Y$}_x\vert^2
=1-\left(\frac{p-x}{\vert x-p\vert}\cdot\mbox{\boldmath $\nu$}_x\right)^2,
$$
it follows from Lemma 4.2 that, for all $x\in\partial D\cap B_{d_{\partial D}(p)+\delta}(p)$
$$
\displaystyle
\vert\mbox{\boldmath $Y$}_x\vert
\le\sqrt{2\eta}.
\tag {4.13}
$$

Since $\mbox{\boldmath $a$}_1$ and $\mbox{\boldmath $a$}_2$ are
linearly independent, it is easy to see that, for all
$x\in\partial D$ we have $\sum_{j=1}^2\vert\mbox{\boldmath
$\nu$}_x\times(\mbox{\boldmath $a$}_j\times\mbox{\boldmath
$\nu$}_x)\vert^2>0$. Since $\partial D$ is compact, this yields
that $\min_{x\in\partial D}\sum_{j=1}^2\vert\mbox{\boldmath
$\nu$}_x\times(\mbox{\boldmath $a$}_j\times\mbox{\boldmath
$\nu$}_x)\vert^2>0$. Choose a sufficiently small $\eta>0$ in such
a way that
$$
\displaystyle
\eta<\min_{x\in\partial D}\sum_{j=1}^2\vert\mbox{\boldmath $\nu$}_x\times(\mbox{\boldmath $a$}_j\times\mbox{\boldmath $\nu$}_x)\vert^2.
$$

We have
$$\begin{array}{l}
\displaystyle
\,\,\,\,\,\,\,
-\left(
P_x\,\mbox{\boldmath $X$}_x
-Q_x\,\mbox{\boldmath $Y$}_x
\right)\cdot\left\{
\left(\mbox{\boldmath $\nu$}_x\cdot\frac{x-p}{\vert x-p\vert}\right)\mbox{\boldmath $X$}_x
-(\mbox{\boldmath $\nu$}_x\cdot\mbox{\boldmath $a$}_j)\mbox{\boldmath $Y$}_x\right\}\\
\\
\displaystyle
=P_x\left(\mbox{\boldmath $\nu$}_x\cdot\frac{p-x}{\vert x-p\vert}\right)\vert\mbox{\boldmath $X$}_x\vert^2
\\
\\
\displaystyle
\,\,\,
-Q_x(\mbox{\boldmath $\nu$}_x\cdot\mbox{\boldmath $a$}_j)\vert\mbox{\boldmath $Y$}_x\vert^2
+\left\{
P_x(\mbox{\boldmath $\nu$}_x\cdot\mbox{\boldmath $a$}_j)-Q_x\left(\mbox{\boldmath $\nu$}_x\cdot\frac{p-x}{\vert x-p\vert}\right)
\right\}
\mbox{\boldmath $Y$}_x\cdot\mbox{\boldmath $X$}_x.
\end{array}
\tag {4.14}
$$

Write
$$\begin{array}{ll}
\displaystyle
P_x & \displaystyle
=\lambda(x)
\left\{1+\frac{1}{\tau\sqrt{\mu\epsilon}}
\left(\frac{1}{\vert x-p\vert}+\frac{1}{\tau\sqrt{\mu\epsilon}\vert x-p\vert^2}\right)\right\}
\\
\\
\displaystyle
& \displaystyle
\,\,\,
-\sqrt{\frac{\epsilon}{\mu}}
\,\left(1+\frac{1}{\tau\sqrt{\mu\epsilon}\vert x-p\vert}\right)\mbox{\boldmath $\nu$}\cdot\frac{p-x}{\vert x-p\vert}
\\
\\
\displaystyle
& \displaystyle
=\left(1+\frac{1}{\tau\sqrt{\mu\epsilon}\vert x-p\vert}\right)
\left(\lambda(x)-\sqrt{\frac{\epsilon}{\mu}}\,\mbox{\boldmath $\nu$}\cdot\frac{p-x}{\vert x-p\vert}\right)
+\frac{\lambda(x)}{(\tau\sqrt{\mu\epsilon})^2\vert x-p\vert^2}.
\end{array}
\tag {4.15}
$$
Since we have, for all $x\in\partial D$
$$\begin{array}{c}
\displaystyle
\lambda(x)-\sqrt{\frac{\epsilon}{\mu}}\,\mbox{\boldmath $\nu$}_x\cdot\frac{p-x}{\vert x-p\vert}
=\left(\lambda(x)-\sqrt{\frac{\epsilon}{\mu}}\,\,\right)
+\sqrt{\frac{\epsilon}{\mu}}
\left(1-\mbox{\boldmath $\nu$}_x\cdot\frac{p-x}{\vert x-p\vert}\right)
\\
\\
\displaystyle
\ge\lambda(x)-\sqrt{\frac{\epsilon}{\mu}}
\ge C,
\end{array}
$$
it follows from these and Lemma 4.2 for $\eta<1/2$ we obtain, for all $\tau>0$
and $x\in\partial D\cap B_{d_{\partial D}(p)+\delta}(p)$
$$\begin{array}{c}
\displaystyle
P_x\left(\mbox{\boldmath $\nu$}_x\cdot\frac{p-x}{\vert x-p\vert}\right)\ge\frac{C}{2}.
\end{array}
\tag {4.16}
$$
Using (4.13) and simply estimating from above, we have
$$\begin{array}{l}
\displaystyle
\,\,\,\,\,\,
\left\vert-Q_x(\mbox{\boldmath $\nu$}_x\cdot\mbox{\boldmath $a$}_j)\vert\mbox{\boldmath $Y$}_x\vert^2
+\left\{
P_x(\mbox{\boldmath $\nu$}_x\cdot\mbox{\boldmath $a$}_j)-Q_x\left(\mbox{\boldmath $\nu$}_x\cdot\frac{p-x}{\vert x-p\vert}\right)
\right\}
\mbox{\boldmath $Y$}_x\cdot\mbox{\boldmath $X$}_x\right\vert
\\
\\
\displaystyle
\le C'(\eta+\sqrt{\eta}\,\vert\mbox{\boldmath $X$}_x\vert).
\end{array}
$$
Applying this together with (4.16) to (4.14), we have, for all $\kappa>0$
$$\begin{array}{l}
\displaystyle
\,\,\,\,\,\,
-\left(
P_x\,\mbox{\boldmath $X$}_x
-Q_x\,\mbox{\boldmath $Y$}_x
\right)\cdot\left\{
\left(\mbox{\boldmath $\nu$}_x\cdot\frac{x-p}{\vert x-p\vert}\right)\mbox{\boldmath $X$}_x
-(\mbox{\boldmath $\nu$}_x\cdot\mbox{\boldmath $a$}_j)\mbox{\boldmath $Y$}_x\right\}\\
\\
\displaystyle
\ge \frac{C}{2}\vert\mbox{\boldmath $X$}_x\vert^2
-C'(\eta+\sqrt{\eta}\,\vert\mbox{\boldmath $X$}_x\vert)\\
\\
\displaystyle
\ge
\left(\frac{C}{2}-\frac{C'\kappa}{2}\right)\,\vert\mbox{\boldmath $X$}_x\vert^2
-C'\left(\eta+\frac{\kappa^{-1}\eta}{2}\right).
\end{array}
$$
Thus, letting $\kappa=C/(2C')$, we obtain, for all $x\in\partial D\cap B_{d_{\partial D}(p)+\delta}(p)$ and $\tau>0$
$$\displaystyle
-\left(
P_x\,\mbox{\boldmath $X$}_x
-Q_x\,\mbox{\boldmath $Y$}_x
\right)\cdot\left\{
\left(\mbox{\boldmath $\nu$}_x\cdot\frac{x-p}{\vert x-p\vert}\right)\mbox{\boldmath $X$}_x
-(\mbox{\boldmath $\nu$}_x\cdot\mbox{\boldmath $a$}_j)\mbox{\boldmath $Y$}_x\right\}
\ge \frac{C}{4}\vert\mbox{\boldmath $X$}_x\vert^2-C_2\eta,
$$
where $C_2=C'(C+C')/C$.

Now it follows from this and (4.11) that, for all $x\in\partial D\cap B_{d_{\partial D}(p)+\delta}(p)$
$$\begin{array}{c}
\displaystyle
\sum_{j=1}^2\frac{1}{c}\left\{c\,\mbox{\boldmath $\nu$}_x\times(\mbox{\boldmath$V$}_{e,\,j}^0(x)\times\mbox{\boldmath $\nu$}_x)
-(\nabla\times\mbox{\boldmath $V$}_{e,\,j}^0(x))\times\mbox{\boldmath $\nu$}_x\right\}
\cdot(\nabla\times\mbox{\boldmath $V$}_{e,\,j}^0(x))\times\mbox{\boldmath $\nu$}_x
\\
\\
\displaystyle
\ge
(\lambda(x))^{-1}\sqrt{\mu\epsilon}\,\tau K(\tau)^2(\tilde{f}(\tau))^2v(x)^2
\left(\frac{C}{4}\sum_{j=1}^2\vert\mbox{\boldmath $\nu$}_x\times(\mbox{\boldmath $a$}_j\times\mbox{\boldmath $\nu$}_x)\vert^2-2C_2\eta
\right).
\end{array}
\tag {4.17}
$$
Now re-choosing $\eta$ in such a way that
$$\displaystyle
\eta<\frac{C}{16C_2}\min_{x\in\partial D}\sum_{j=1}^2
\vert\mbox{\boldmath $\nu$}_x\times(\mbox{\boldmath $a$}_j\times\mbox{\boldmath $\nu$}_x)\vert^2,
$$
we have
$$\begin{array}{c}
\displaystyle
(\lambda(x))^{-1}\sqrt{\mu\epsilon}\,\tau K(\tau)^2(\tilde{f}(\tau))^2v(x)^2
\left(\frac{C}{4}\sum_{j=1}^2\vert\mbox{\boldmath $\nu$}_x\times(\mbox{\boldmath $a$}_j\times\mbox{\boldmath $\nu$}_x)\vert^2-2C_4\eta\right)\\
\\
\displaystyle
\ge
(\Vert \lambda\Vert_{L^{\infty}((\partial D)_{\delta}(p)})^{-1}\sqrt{\mu\epsilon}\,\tau K(\tau)^2(\tilde{f}(\tau))^2v(x)^2
\frac{C}{8}\min_{x\in\partial D}\sum_{j=1}^2\vert\mbox{\boldmath $\nu$}_x\times(\mbox{\boldmath $a$}_j\times\mbox{\boldmath $\nu$}_x)\vert^2.
\end{array}
$$
Therefore this together with (4.17) yields
$$\begin{array}{l}
\displaystyle
\,\,\,\,\,\,
\sum_{j=1}^2\frac{1}{\mu\epsilon}
\int_{\partial D\cap B_{d_{\partial D}(p)+\delta}(p)}
\frac{1}{c}\left\{c\,\mbox{\boldmath $\nu$}\times(\mbox{\boldmath$V$}_{e,\,j}^0\times\mbox{\boldmath $\nu$})
-(\nabla\times\mbox{\boldmath $V$}_{e,\,j}^0)\times\mbox{\boldmath $\nu$}\right\}
\cdot(\nabla\times\mbox{\boldmath $V$}_{e,\,j}^0)\times\mbox{\boldmath $\nu$}\,dS
\\
\\
\displaystyle
\ge C_3(\Vert\lambda\Vert_{L^{\infty}(\partial D)}\sqrt{\mu\epsilon}\,)^{-1}\tau K(\tau)^2(\tilde{f}(\tau))^2
\int_{\partial D\cap B_{d_{\partial D}(p)+\delta}(p)}v(x)^2dS,
\end{array}
\tag {4.18}
$$
where
$$\displaystyle
C_3=\frac{C}{8}
\min_{x\in\partial D}\sum_{j=1}^2\vert\mbox{\boldmath $\nu$}_x\times(\mbox{\boldmath $a$}_j\times\mbox{\boldmath $\nu$}_x)\vert^2.
$$
Using the same argument done in the proof of Lemma 2.2 in \cite{IE2}, we know that
there exists a positive constant  $C''$ such that, for all $\tau>>1$
$$\displaystyle
\int_{\partial D\cap B_{d_{\partial D}(p)+\delta}(p)}v(x)^2dS
\ge C''\tau^{-4}e^{-2\tau \sqrt{\mu\epsilon}d_{\partial D}(p)}.
$$
And it is easy to see that, as $\tau\longrightarrow\infty$
$$\displaystyle
K(\tau)\sim\tau^{-1}\frac{\displaystyle
\eta e^{\tau\eta\sqrt{\mu\epsilon}}}
{2\epsilon}.
\tag {4.19}
$$
Therefore we have
$$\displaystyle
\tau K(\tau)^2\int_{\partial D\cap B_{d_{\partial D}(p)+\delta}(p)}v(x)^2dS
\ge C'''\tau^{-5}e^{-2\tau\sqrt{\mu\epsilon}\,\text{dist}\,(D,B)},
\tag {4.20}
$$
where $C'''$ is a positive constant and $\tau>>1$. Since
$\mbox{\boldmath $V$}_e^0$ together with its derivatives on
$\partial D\setminus B_{d_{\partial D}(p)+\delta}(p)$ is decaying
as $e^{-\tau\sqrt{\mu\epsilon}\text{dist}\,(D,B)}
e^{-\tau\sqrt{\mu\epsilon}\delta}$ (see (4.10) and note that
$d_{\partial D}(p)-\eta=\text{dist}\,(D,B)$), (4.1), (4.18),
(4.20) and (1.11) yields (i) with $\rho=2\gamma+5$.

\subsection{Proof of (ii)}

Write
$$\begin{array}{l}
\displaystyle
\,\,\,\,\,\,
\left(
P_x\,\mbox{\boldmath $X$}_x
-Q_x\,\mbox{\boldmath $Y$}_x
\right)
\cdot
\left\{
A\,\mbox{\boldmath $X$}_x
-B\,
\left(\mbox{\boldmath $a$}\cdot\frac{x-p}{\vert x-p\vert}\right)
\mbox{\boldmath $Y$}_x
\right\}\\
\\
\displaystyle
=P_xA\vert\mbox{\boldmath $X$}\vert^2
\\
\\
\displaystyle
\,\,\,
-\left\{Q_xA+P_xB\left(\mbox{\boldmath $a$}\cdot\frac{x-p}{\vert x-p\vert}\right)\right\}\mbox{\boldmath $X$}_x\cdot\mbox{\boldmath $Y$}_x
+Q_xB\left(\mbox{\boldmath $a$}\cdot\frac{x-p}{\vert x-p\vert}\right)\vert\mbox{\boldmath $Y$}_x\vert^2.
\end{array}
\tag {4.21}
$$
It follows from Lemma 4.2 that, for all $x\in\partial D\cap B_{d_{\partial D}(p)+\delta}(p)$
$$\begin{array}{ll}
\displaystyle
\lambda(x)-\sqrt{\frac{\epsilon}{\mu}}\,\mbox{\boldmath $\nu$}_x\cdot\frac{p-x}{\vert x-p\vert}
&
\displaystyle
=\left(\lambda(x)-\sqrt{\frac{\epsilon}{\mu}}\,\,\right)
+\sqrt{\frac{\epsilon}{\mu}}
\left(1-\mbox{\boldmath $\nu$}_x\cdot\frac{p-x}{\vert x-p\vert}\right)
\\
\\
\displaystyle
&
\displaystyle
\le -C+\sqrt{\frac{\epsilon}{\mu}}\,\eta.
\end{array}
$$
Thus, choosing $\eta$ in such a way that
$$\displaystyle
\eta<\frac{C}{2}\sqrt{\frac{\mu}{\epsilon}},
$$
we have
$$\displaystyle
\lambda(x)-\sqrt{\frac{\epsilon}{\mu}}\,\mbox{\boldmath $\nu$}_x\cdot\frac{p-x}{\vert x-p\vert}
=\left(\lambda(x)-\sqrt{\frac{\epsilon}{\mu}}\,\,\right)
+\sqrt{\frac{\epsilon}{\mu}}
\left(1-\mbox{\boldmath $\nu$}_x\cdot\frac{p-x}{\vert x-p\vert}\right)
\le -\frac{C}{2}.
$$
Therefore from (4.15) we obtain
$$\displaystyle
P_x\le -\frac{C}{2}+\frac{\Vert\lambda\Vert_{L^{\infty}(\partial D)}}{(\tau\sqrt{\mu\epsilon})^2(d_{\partial D}(p))^2}.
$$
Thus, choosing $\tau_0>0$ in such a way that
$$\displaystyle
\frac{\Vert\lambda\Vert_{L^{\infty}(\partial D)}}{(\tau_0\sqrt{\mu\epsilon})^2(d_{\partial D}(p))^2}=\frac{C}{4},
$$
we have, for all $\tau\ge\tau_0$ and all $x\in\partial D\cap B_{d_{\partial D}(p)+\delta}(p)$
$$
\displaystyle
P_x\le -\frac{C}{4}.
$$
Since $A\ge 1$, this yields
$$\displaystyle
P_xA\le-\frac{C}{4}.
\tag {4.22}
$$
It is easy to see that we have
$$\begin{array}{l}
\displaystyle
\,\,\,\,\,\,
\left\vert-\left\{Q_xA+P_xB\left(\mbox{\boldmath $a$}\cdot\frac{x-p}{\vert x-p\vert}\right)\right\}\mbox{\boldmath $X$}_x\cdot\mbox{\boldmath $Y$}_x
+Q_xB\left(\mbox{\boldmath $a$}\cdot\frac{x-p}{\vert x-p\vert}\right)\vert\mbox{\boldmath $Y$}_x\vert^2\right\vert
\\
\\
\displaystyle
\le
C_4(\eta+\sqrt{\eta}\vert\mbox{\boldmath $X$}_x\vert).
\end{array}
$$
Applying this together with (4.22) to (4.21), we obtain
$$\begin{array}{l}
\displaystyle
\left(
P_x\,\mbox{\boldmath $X$}_x
-Q_x\,\mbox{\boldmath $Y$}_x
\right)
\cdot
\left\{
A\,\mbox{\boldmath $X$}_x
-B\,
\left(\mbox{\boldmath $a$}\cdot\frac{x-p}{\vert x-p\vert}\right)
\mbox{\boldmath $Y$}_x
\right\}
\\
\\
\displaystyle
\le -\frac{C}{4}\vert\mbox{\boldmath $X$}_x\vert^2+C_4(\eta+\sqrt{\eta}\vert\mbox{\boldmath $X$}_x\vert)\\
\\
\displaystyle
\le
-\left(\frac{C}{4}-\frac{C_4\kappa}{2}\right)\vert\mbox{\boldmath $X$}_x\vert^2
+C_4\eta\left(1+\frac{1}{2\kappa}\right),
\end{array}
$$
where $\kappa$ is an arbitrary positive number.  Thus letting $\kappa=C/(4C_4)$, we obtain
$$\begin{array}{c}
\displaystyle
\left(
P_x\,\mbox{\boldmath $X$}_x
-Q_x\,\mbox{\boldmath $Y$}_x
\right)
\cdot
\left\{
A\,\mbox{\boldmath $X$}_x
-B\,
\left(\mbox{\boldmath $a$}\cdot\frac{x-p}{\vert x-p\vert}\right)
\mbox{\boldmath $Y$}_x
\right\}
\le -\left(\frac{C}{8}\vert\mbox{\boldmath $X$}_x\vert^2-C_5\eta\right),
\end{array}
$$
where $C_5=C_4(C+2C_4)/C$.

Therefore this together with (4.12) yields that, for all $x\in\partial D\cap B_{d_{\partial D}(p)+\delta}(p)$
and all $\tau\ge\tau_0$
$$\begin{array}{l}
\displaystyle
\,\,\,\,\,\,
\sum_{j=1}^2\left\{c\,\mbox{\boldmath $\nu$}_x\times(\mbox{\boldmath$V$}_{e,\,j}^0(x)\times\mbox{\boldmath $\nu$}_x)
-(\nabla\times\mbox{\boldmath $V$}_{e,\,j}^0(x))\times\mbox{\boldmath $\nu$}_x\right\}
\cdot\mbox{\boldmath $\nu$}_x\times(\mbox{\boldmath $V$}_{e,\,j}^0(x)\times\mbox{\boldmath $\nu$}_x)
\\
\\
\displaystyle
\le
-\mu\tau K(\tau)^2(\tilde{f}(\tau))^2v(x)^2
\left(\frac{C}{8}\sum_{j=1}^2\vert\mbox{\boldmath $\nu$}_x\times(\mbox{\boldmath $a$}_j\times\mbox{\boldmath $\nu$}_x)\vert^2-2C_5
\eta\right).
\end{array}
\tag {4.23}
$$
Now re-choosing $\eta$ in such a way that
$$\displaystyle
\eta<\frac{C}{32C_5}\min_{x\in\partial D}\sum_{j=1}^2
\vert\mbox{\boldmath $\nu$}_x\times(\mbox{\boldmath $a$}_j\times\mbox{\boldmath $\nu$}_x)\vert^2,
$$
from (4.23) we obtain
$$
\begin{array}{l}
\displaystyle
\,\,\,\,\,\,
\sum_{j=1}^2\frac{1}{\mu\epsilon}
\int_{\partial D\cap B_{d_{\partial D}(p)+\delta}(p)}
\left\{c\,\mbox{\boldmath $\nu$}\times(\mbox{\boldmath$V$}_{e,\,j}^0\times\mbox{\boldmath $\nu$})
-(\nabla\times\mbox{\boldmath $V$}_{e,\,j}^0)\times\mbox{\boldmath $\nu$}\right\}
\cdot\mbox{\boldmath $\nu$}\times(\mbox{\boldmath $V$}_{e,\,j}^0\times\mbox{\boldmath $\nu$})\,dS
\\
\\
\displaystyle
\le
-C_6\epsilon^{-1}\tau K(\tau)^2(\tilde{f}(\tau))^2
\int_{\partial D\cap B_{d_{\partial D}(p)+\delta}(p)}v(x)^2dS,
\end{array}
$$
where
$$\displaystyle
C_6=\frac{C}{16}
\min_{x\in\partial D}\sum_{j=1}^2\vert\mbox{\boldmath $\nu$}_x\times(\mbox{\boldmath $a$}_j\times\mbox{\boldmath $\nu$}_x)\vert^2.
$$
Hereafter the procedure is completely same as that of the proof of (i).

\subsection{Proof of Lemma 2.3}

From (4.1), (4.11) and (4.19) we have, as $\tau\longrightarrow\infty$
$$\displaystyle
\tilde{J}_e(\tau)
=O\left(\tau^{-1}(\tilde{f}(\tau))^2e^{2\tau\sqrt{\mu\epsilon}}\int_{\partial D}v^2dS\right)+O(\tau^{-1/2}e^{-\tau T})
$$
and also (4.2), (4.11) and (4.19) give
$$\displaystyle
\tilde{J}_e(\tau)
+\frac{1}{\mu\epsilon}\int_{\partial D}c\vert\mbox{\boldmath $V$}_{em}\vert^2\,dS
=O\left(\tau^{-1}(\tilde{f}(\tau))^2e^{2\tau\sqrt{\mu\epsilon}}\int_{\partial D}v^2dS\right)+O(\tau^{-1/2}e^{-\tau T}).
$$
Then, it follows from these, (2.1) and (2.2) that (2.5) is valid.

\section{Conclusions and further problems}

We have succeeded in extending the previous result in \cite{IE2}
for the scalar wave equation with the dissipative boundary
condition to the Maxwell system with the Leontovich boundary
condition.   It can be also considered as an extension of Theorem
1.1 in  \cite{IEM} for the Maxwell system with the perfect
conductivity condition to the Leontovich boundary condition.  The
main difference is to introduce an indicator function which
employs two sets of electric fields corresponding to input sources
oriented to two independent directions:
$$\displaystyle
\tau\longmapsto\sum_{j=1}^2 I_{\mbox{\boldmath $f$}_j}(\tau).
\tag {5.1}
$$

Technically the previous indicator function of Theorem 1.1 in \cite{IEM}
is based on the behaviour of a {\it volume integral} of $\mbox{\boldmath $V$}_e^0$ over the obstacle
which is reduced to that of the fundamental solution of the modified Helmholz equation.
In contrast to the previous case, as a consequence of the presence of the impedance $\lambda$
on the obstacle, the asymptotic behaviour of each of two components of the new indicator function
is governed by the  behaviour of  a {\it surface integral} of functions involving $\mbox{\boldmath $V$}_{e}^0$
together with derivatives.  Since $\mbox{\boldmath $V$}_e^0$ has a {\it directivity}, it will be difficult to
obtain the same result as Theorem 1.1 in \cite{IEM} by using only a single component of the indicator function
without assuming some restrictive condition for the source direction like (9) in \cite{IEM}.

Our finding is: to obtain the distance of a given point to an unknown obstacle
together with a qualitative state of the surface of the obstacle
{\it without} any restriction on the orientation of the source
it is enough to make use of two sets of electric fields generated by
two linearly independent input sources supported on a common open ball.

In future research \cite{IC} we will consider: extract information
about the {\it curvatures} and impedance at the points on the
surface of the obstacle nearest to the center of the support of
input sources from the asymptotic behaviour of indicator function
(5.1).   This can be considered as an extension of Theorem 1.2 in
\cite{IEM}. And also it would be interested to consider the
extension of the results in \cite{IW} to the Maxwell system. See
\cite{IR} for a survey on other results using the enclosure method
in the time domain together with other open problems.

$$\quad$$

\centerline{{\bf Acknowledgment}}

The author was partially supported by Grant-in-Aid for
Scientific Research (C)(No. 25400155) of Japan  Society for
the Promotion of Science.

$$\quad$$

\vskip1cm
\noindent
e-mail address

ikehata@amath.hiroshima-u.ac.jp

\end{document}